\documentclass[12pt]{article}%
\usepackage{amsmath}
\usepackage{amssymb}
\usepackage{dsfont}
\usepackage{cite}
\newsymbol \blackbox 1004
\newcommand{\eh}{\hfill}\newlength{\sperr}

\newenvironment{proof}{{\settowidth{\sperr}{\bf\rm
Proof}%
\par\addvspace{0.3cm}\noindent\parbox[t]{1.3\sperr}
{\bf\rm P\eh r\eh o\eh o\eh f\eh }%
}}{\nopagebreak\mbox{}
$\blackbox$\par\addvspace{0.3cm}}

\def\nn{\nonumber}

\def\a{\alpha}
\def\b{\beta}

\def\de{\delta}

\def\Lam{\Lambda}

\def\la{\lambda}
\def\om{\omega}

\def\Th{\Theta}
\def\ze{\zeta}
\def\Up{\Upsilon}
\def\up{\upsilon}
\def\vp{\varphi}
\def\vt{\vartheta}
\def\ve{\varepsilon}
\def\wh{\widehat}
\def\wt{\widetilde}
\def\ov{\overline}
\def\br{\breve}

\def\p{\partial}

\def\BC{{\mathbb C}}

\def\BR{{\mathbb R}}
\def\BN{{\mathbb N}}
\def\clp{{\mathcal P}}
\def\cla{{\mathcal A}}

\def\cle{{\mathcal E}}
\def\clf{{\mathcal F}}
\def\clg{{\mathcal G}}
\def\clh{{\mathcal H}}

\def\clk{{\mathcal K}}

\def\cln{{\mathcal N}}
\def\clu{{\mathcal U}}
\def\clw{{\mathcal W}}
\def\clq{{\mathcal Q}}
\def\clr{\mathcal{R}}

\def\clv{{\mathcal V}}

\def\clz{{\mathcal Z}}

\newcommand{\E}{\mathrm{e}}
\newcommand{\I}{\mathrm{i}}
\def\mf{\mathfrak}

\def\bB{{\bf B}}

\newtheorem{Pa}{Paper}[section]
\newtheorem{Tm}[Pa]{{\bf Theorem}}
\newtheorem{La}[Pa]{{\bf Lemma}}
\newtheorem{Cy}[Pa]{{\bf Corollary}}
\newtheorem{Rk}[Pa]{{\bf Remark}}

\newtheorem{Ee}[Pa]{{\bf Example}}
\newtheorem{Dn}[Pa]{{\bf Definition}}
\newtheorem{Nn}[Pa]{{\bf Notation}}
\newtheorem{Pn}[Pa]{{\bf Proposition}}

\title{Generalised canonical systems related to
matrix string equations: corresponding structured operators and high-energy asymptotics of  the Weyl functions}

\author{Alexander Sakhnovich}

\date{}
\begin{document}
\maketitle

\begin{abstract}  We obtain high energy asymptotics of  Titchmarsh-Weyl functions of the generalised
canonical systems generalising in this way a seminal Gesztesy-Simon result. The 
matrix valued analog of the amplitude function satisfies in this case an interesting  new identity.
The corresponding structured operators are studied as well. Application to a procedure of
solving an important inverse problem is presented in Appendix.
\end{abstract}


\vspace{0.2em}


\section{Introduction} \label{intro}
\setcounter{equation}{0}
Canonical systems have the form
\begin{align} &       \label{1.1}
w^{\prime}(x,\la)=\I \la J H(x)w(x,\la), \quad J:=\begin{bmatrix} 0 & I_p \\ I_p & 0\end{bmatrix} \quad \Big(w^{\prime}:=\frac{d}{dx}w\Big),
 \end{align} 
where $\I$ is the  imaginary unit ($\I^2=-1$), $\la$ is the so called spectral parameter,
$I_{p}$ is the $p \times p$ $(p\in \BN)$ identity
matrix,  $\BN$ stands for the set of positive integer numbers,  $H(x)$ is a $2p \times 2p$  matrix function (matrix valued function),
and 
$H(x) \geq 0$ (that is, the matrices $H(x)$ are self-adjoint and the eigenvalues of $H(x)$ are nonnegative).
Canonical systems are important objects of analysis, being perhaps the most important class of  the one-dimensional
Hamiltonian systems and including (as subclasses) several classical equations. They have been actively studied in many
already classical as well as in various recent works (see, e.g., \cite{ArD, dBr,  EKT,  GoKr, GolMi, Langer, Mog,  Rom, RW, Rov,  ALS17, ALSstring,
SaSaR, SaL2-, LA94, SaL2, Su, Wor} 
and numerous references
therein). 

In most works on canonical systems,  a somewhat simpler  case of $2 \times 2$ Hamiltonians $H(x)$ (i.e, the case $p=1$) is dealt with.
In particular, the trace normalisation $\mathrm{tr} \, H(x)\equiv 1$ may be successfully used in the case of $p=1$.
The cases with other values of $p$ ($p>1$) are equally important but more complicated and less studied. 
\begin{Rk}\label{RkRem} The fundamental results by L. de Branges and by M.G. Krein on canonical systems,
where the Hamiltonians are $2\times 2$  trace-normalised matrix functions with real-valued entries,
are well presented $($and usefully reformulated sometimes$)$ in the recent book \cite{Rem}.
In particular, the basic de Branges result on the one to one correspondence between generalized Herglotz  functions and  the above-mentioned canonical systems
is presented there.  An interesting
research on the absolutely continuous spectrum  by the author $($C.~Remling$)$   himself  is also contained in \cite{Rem}.
However, the procedure for solving inverse problem is missing in this excellent book.\end{Rk}

In our paper, we consider  
{\it generalised canonical systems}
 \begin{equation}        \label{1.2}
 w^{\prime}(x,\la)=\I \la j H(x)w(x,\la), \,\, H(x)\geq 0, \,\, j:=\begin{bmatrix} I_{m_1} & 0 \\  0 & -I_{m_2}\end{bmatrix} \,\, (m_1+m_2=:m),
 \end{equation} 
where  
$x \geq 0$, $\, m_1,m_2 \in \BN$, 
 and  $H$ is an $m\times m$ locally integrable matrix function. More precisely, we consider mostly the case 
 of the Hamiltonians  $H(x)$ of the form
 \begin{align}&\label{1.3}
 H(x)=\b(x)^*\b(x),
 \end{align}
where $\b$ are $m_2 \times m$ matrix functions and $\b(x)^*$ is the matrix adjoint to $\b(x)$.  
Systems \eqref{1.2}, \eqref{1.3} are studied on $[0,r]$
or $[0,\infty)$.
We assume that $\b(x)\in \clu^{m_2\times m}[0,r]$, where
 \begin{align}&\label{1.3+}
\clu^{p\times q}[0,r]=\big\{\clg: \,\, \clg^{\prime}(x)\equiv \clg^{\prime}(0)+\int_0^x \clg^{\prime\prime}(t)dt, \quad
\clg^{\prime\prime}\in L_2^{p\times q}(0,r)\big\},
 \end{align}
$L_2^{p\times q}(0,r)$ stands for the class of  $p\times q$
matrix functions with square integrable entries (i.e. the entries from $L_2(0, r)$) and $\clg^{\prime}$ is the standard derivative of $\clg$.
We say that $\clg$ in \eqref{1.3+} is two times differentiable
and that 
$ \clg^{\prime\prime}$ satisfying $ \clg^{\prime}(x)\equiv \clg^{\prime}(0)+\int_0^x \clg^{\prime\prime}(t)dt$
is the second derivative of $\clg$. 
\begin{Rk} \label{RkSob} 
Clearly, $\clu^{p\times q}[0,r]$ is the class of $p\times q$ matrix functions, the entries of which belong to the Sobolev class
$W^{2,2}$ on $[0,r]$ $($although Sobolev norms are not used in this work$).$
\end{Rk}
When 
 \begin{align}& \label{1.4-}
\b(x)\in \clu^{m_2\times m}[0,r] 
\end{align}
${\mathrm{for \, all}} \,\, r>0$, we write $\b(x)\in \clu^{m_2\times m}[0,\infty)$.
We  also assume that
 \begin{align}&\label{1.4}
 \b(x)j\b(x)^*\equiv 0, \quad  \b^{\prime}(x)j\b(x)^* \equiv \I I_{m_2}.
 \end{align}

When $m_1=m_2=p$, system \eqref{1.2} is equivalent to the canonical system \eqref{1.1} with a slightly different
Hamiltonian (see Section \ref{WH}).  Moreover, system \eqref{1.2}--\eqref{1.4}
(in this case and under some minor additional conditions) may be transformed into the matrix string equation \cite[Appendix B]{ALSstring}.
The case of matrix Schr\"odinger equations is included in this class (see \cite[\S 2 of Appendix B]{ALSstring}).
Systems \eqref{1.2}--\eqref{1.4} present a class of canonical systems complimentary to the canonical systems corresponding to
Dirac systems. (On the canonical  systems corresponding to
Dirac systems see, e.g., \cite{ALS15, SaSaR} and references therein.)

This paper is an important continuation and development of our article \cite{ALSstring}. Here,
we study the high  energy asymptotics of  Weyl (Titchmarsh--Weyl) functions of  the systems \eqref{1.2}--\eqref{1.4}.
Asymptotics of the Weyl functions is an important topic with interesting results obtained, in particular,
in \cite{Atk, Ev}.  Some fundamental
results  on Weyl functions followed (for the scalar Schr\"odinger equation) in the seminal papers \cite{Si0} by B. Simon and  \cite{GeSi}
by F. Gesztesy and B. Simon (see further discussion in \cite{GeB}). 
Closely related results in terms of spectral functions have been stated in the pioneering note  \cite{Krein0} by M.G.~Krein (unfortunately
without proofs). 
It was later shown by H. Langer \cite{Langer} that the assertions from \cite{Krein0} yield high energy
asymptotics for Weyl functions as well  (at least for the case treated in \cite{Si0}). The corresponding assertions
from \cite{Krein0}  (for the special case of the orthogonal  spectral functions) were also proved in \cite{Langer}.
Similar relations for Weyl functions corresponding to string equations one can find in \cite{ALS88} (see formula
(45) and Statement 7 there based on the work \cite{SaL2-} on $S$-colligations).
The case of canonical systems, such that the Hamiltonians $H(x)$ are $2\times 2$ matrix functions with real-valued
entries and tr $H(x)\equiv 1$, was treated in the interesting papers \cite{KLang, Langer, LW}.  

Related results one can find in \cite{Krein} and in \cite[(32)]{ALS88} for self-adjoint Dirac systems,
and  in \cite[p. 319]{ALS90} for skew-self-adjoint Dirac systems, with various important  developments in
\cite{ClGe, ClGe2,  ALS02}.  The inverse approach to Dirac  systems 
\begin{align} &      \nonumber
y^{\prime}(x, \la )=i (\la j + j V(x))y(x,\la ), \quad V= \begin{bmatrix}
0 & v \\ v^{*} & 0 \end{bmatrix},
\end{align} 
(with $m_1 \times m_2$ potentials $v$)
based on the A-function concept was studied in \cite{GeS}.

The main result of the present paper is Theorem \ref{TmHEA}. This theorem generalises 
(for the case of the canonical systems \eqref{1.2}--\eqref{1.4}) important Theorem~1.1 from \cite{GeSi},
which constitutes, as also mentioned by the authors of \cite{GeSi}, one of the main results
of their seminal work. The result is new even for the $2\times 2$ Hamiltonians (in particular, we do not
require that the entries of $H$ are real-valued). An interesting new identity \eqref{S39'} for the analog of the so called
$A$-amplitude appears here.  We also give an analogue of  Theorem \ref{TmHEA} for the case of $m_1=m_2$ and 
Weyl matrix functions belonging to Herglotz class (see Theorem \ref{TmHEAH}).
\begin{Rk} \label{Rk}  The results on the so called high energy asymptotics of Weyl functions $($of the \eqref{M9} and \eqref{B22} type$)$
have various applications including applications to the local uniqueness and other problems in the inverse spectral theory  \cite{Si0, GeSi, GeSi2}
$($see also \cite{GeS, Langer, Ryb, ALS88} and references therein$)$. The work on the applications of  \eqref{M9} and \eqref{B22}
is in progress. However, the uniqueness and procedure of solving  inverse problem for the important case
$m_1=m_2$  is already  derived from the results of this paper \cite{ALSinvpr}.
\end{Rk}

The next Section \ref{Prel} is called ``Preliminaries". We generalise their some
basic results from \cite{ALSstring} and present an interesting Example \ref{Ee2.3}. Section \ref{Fund} is dedicated to the fundamental
solutions of the canonical systems.
The corresponding operator identities and structured operators are
also of independent interest. They are studied in Section \ref{Struct}.
The obtained results are summed up in Section \ref{Direct}
in the form of Theorem~\ref{TmHEA}. An analogue of Theorem \ref{TmHEA} for the case
$m_1=m_2$ and Weyl matrix functions belonging to Herglotz class as well
as two useful examples are given in Section \ref{WH}.

Some linear similarity
problems are discussed in Appendix \ref{Simil}.

{\it Notations.} Some notations were already introduced in the introduction above.
As usual, $\BR$ stands for the real axis,  $\BC$ stands for the complex plane,
the open upper half-plane is denoted
by $\BC_+$, and  $\ov{a}$ means the complex conjugate of $a$. The notation $\Re(a)$ stands for the real part of $a$,
and $\Im(a)$ denotes the imaginary part of $a$. For a matrix $K$, $\ov{K}$ is the matrix such that all its entries are 
complex conjugates of the corresponding entries of $K$. For a matrix or operator $V$, $V^*$ stands for the adjoint
matrix or operator.

We set $L_2^{p\times1}=L_2^p$, $L_2^{1}=L_2$ and  $\clu^{p\times 1}=\clu^p$. ($L_2^p(0,r)$ and $L_2(0,r)$ stand also for the
corresponding Hilbert spaces of square summable functions.)
The notation $I$ stands for the identity operator. The norm $\|A\|$ of the $n\times n $ matrix $A$
means the norm of $A$ acting in the space $\ell_2^n$ of the sequences of length $n$. 
The class of bounded operators acting from the Hilbert space $\clh_1$ into Hilbert space $\clh_2$ is denoted by
$\bB(\clh_1,\clh_2)$, and we set $\bB(\clh):=\bB(\clh,\clh)$.
 If $R\in \bB\big(L_2^p(0,r)\big)$ and $\Phi(x)\in L_2^{p\times q}(0,r)$, then $R\Phi(x) \in L_2^{p\times q}(0,r)$, that is,
 the operators ($R$ here) are applied to matrix functions columnwise.

\section{Preliminaries} \label{Prel}
\setcounter{equation}{0}
This article may be considered as an important development of our paper \cite{ALSstring}, where
fundamental solutions of the class of canonical systems corresponding to matrix string equations were
studied. Simple generalisations of some  basic results, which were obtained  in \cite{ALSstring} for the case 
$m_1=m_2=p$, are presented in this section.

{\bf 1.} Weyl functions of systems \eqref{1.2} (where $m_1=m_2=p$) were considered in
\cite[Appendix A]{ALSstring}. The definitions, results and all proofs in  \cite[Appendix~A]{ALSstring} remain valid after we switch
from the case $m_1=m_2=p$ to the general case of $j$ given in \eqref{1.2}. (Note that we also switch from $I_p$ to $I_{m_1}$.)
More precisely, we have the following relations,
definitions and statements.

The $m\times m$ fundamental solution $W(x,\la)$ of \eqref{1.2} is normalised  by the condition
\begin{align}& \label{P1}
W(0,\la)=I_m.
\end{align}
It is easy to see that
\begin{align}& \label{P2}
\frac{d}{dx}\big(W(x,\ov{\mu})^*jW(x,\la)\big)=\I({\la - \mu})W(x,\ov{\mu})^*H(x)W(x,\la).
\end{align}
In particular, we have the equalities
\begin{align}& \label{P3}
W(r,\ov{\la})^*jW(r,\la)\equiv  j\equiv W(r,\la)jW(r,\ov{\la})^*,
\end{align}
and we set
\begin{align} & \label{P4}
\clw(r,\la)=\{\clw_{ik}(r,\la)\}_{i,k=1}^2:=jW(r,\ov{\la})^*j=W(r,\la)^{-1} \quad (r\geq 0),
\end{align}
where the blocks $\clw_{ik}$ have the same dimensions as the corresponding blocks of $j$ in \eqref{1.2}.

Pairs of meromorphic in $\BC_+$, $m_k\times m_1$ matrix functions $\clp_k(\la)$ $(k=1,2)$ such that
\begin{align} & \label{P6}
\clp_1(\la)^*\clp_1(\la)+\clp_2(\la)^*\clp_2(\la)>0, \quad  \begin{bmatrix}\clp_1(\la)^* & \clp_2(\la)^*\end{bmatrix}j \begin{bmatrix}\clp_1(\la) \\ \clp_2(\la)\end{bmatrix}\geq 0
\end{align}
(where the first inequality holds in one point (at least) of $\BC_+$ and the second inequality holds in all the points of analyticity
of $\clp_1$ and $\clp_2$),  
are called
{\it nonsingular, with property-$j$}.
\begin{Nn} The notation $\cln(r)$ stands for the set of matrix functions of the form
\begin{align} \nn
\phi(r,\la)=&\big(\clw_{21}(r,\la)\clp_1(\la)+\clw_{22}(r,\la)\clp_2(\la)\big)
\\ & \label{P7}
\times\big(\clw_{11}(r,\la)\clp_1(\la)+\clw_{12}(r,\la)\clp_2(\la)\big)^{-1},
\end{align}
where the pairs $\{\clp_1,\clp_2\}$ are nonsingular, with property-$j$.
\end{Nn}

If $\phi(\la)\in \cln(r)$ we have
\begin{align} & \label{P7+}
\begin{bmatrix}I_{m_1} & \phi(\la)^*\end{bmatrix}\mf{A}(r,\la)\begin{bmatrix}I_{m_1} \\ \phi(\la)\end{bmatrix}\geq 0 \quad {\mathrm{for}} \quad
\mf{A}(r,\la):=W(r,\la)^*jW(r,\la).
\end{align}

The matrix functions $\clw_{11}(r,\la)\clp_1(\la)+\clw_{12}(r,\la)\clp_2(\la)$ in \eqref{P7} are invertible (excluding, possibly, isolated
points  $\la \in \BC_+$), and the functions $\phi(r,\la)$ are holomorphic and contractive in $\BC_+$. That is, only removable
singularities of $\phi(r,\la)$ are possible in $\BC_+$.

\begin{Dn}\label{DnW}  Matrix functions $\phi(\la)\in \cln(r)$
are called Weyl $($Titchmarsh--Weyl$)$ functions of the  generalised canonical system \eqref{1.2} on $[0,r]$, where $0<r<\infty$.
Matrix functions $\vp(\la)$ such that
\begin{align} & \label{P8}
\vp(\la)\in \bigcap_{r>0}\cln(r)
\end{align}
are called Weyl   functions of the system  \eqref{1.2} on $[0, \, \infty)$.
\end{Dn}
\begin{Ee}\label{Ee2.3} Consider canonical system for the case
\begin{align} & \label{E10}
\b(x)=\begin{bmatrix}\E^{\I c x}\a & \E^{-\I c x}I_{m_2}\end{bmatrix}, \quad \a \a^*=I_{m_2}\,\, (m_1\geq m_2), \quad c=1/2 ,
\end{align}
where $\a$ is an $m_2\times m_1$ matrix. Clearly, \eqref{1.4} is valid for this $\b$. It follows from \eqref{E10} that
\begin{align} & \label{E11}
H(x)=\b(x)^*\b(x)=\E^{-\I c xj}\clk \E^{\I c xj}, \quad \clk:=\begin{bmatrix}\a^*\a & \a^* \\ \a & I_{m_2}\end{bmatrix}.
\end{align}
One easily checks directly that $W(x,\la)=\E^{-\I c xj}\E^{\I  x(\la j \clk+cj)}$, that is,
\begin{align} & \label{E12}
\clw(x,\la)=\E^{-\I  x(\la j \clk+cj)}\E^{\I c xj}.
\end{align}
Further in this example, we assume that
\begin{align} & \label{E13}
m_1=m_2=p.
\end{align}
It is easily checked that under assumptions \eqref{E10} and \eqref{E13} we have
\begin{align} & \label{E14}
\la j \clk +c j=\Lam\begin{bmatrix}\ze_1 I_{p} & 0 \\0 & \ze_2 I_p\end{bmatrix}\Lam^{-1}, \quad \Lam:=\begin{bmatrix}\a^* & \a^* \\ \ze_1 I_{p} & \ze_2 I_p \end{bmatrix},
\end{align}
where $\ze_1$ and $\ze_2$ are the roots of the quadratic equation
\begin{align}& \label{E14+}
\ze^2+\Big(\frac{1}{\la}+2\Big)\ze+1=0, \quad {\mathrm{that\,\, is}},
\\ & \label{E15}
\ze_1(\la)=-1-\frac{1}{2\la}+\frac{1}{\la}\sqrt{\la+\frac{1}{4}}, \quad \ze_2(\la)=-1-\frac{1}{2\la}-\frac{1}{\la}\sqrt{\la+\frac{1}{4}},
\end{align}
and the branch of the square root in \eqref{E15} belongs the quadrant $\Re( z)>0$,  $\, \Im(z)>0$. Relations
\eqref{E12} and \eqref{E14} yield:
\begin{align} & \label{E16}
\clw(r,\la)=\Lam\begin{bmatrix}\E^{-\I \ze_1 r} I_{p} & 0 \\0 & \E^{-\I \ze_2 r}  I_p\end{bmatrix}\Lam^{-1}\E^{\I c rj}.
\end{align}

Note that in the case \eqref{E13}, the considerations $($of this example$)$ above are similar to the considerations
in \cite[Section 5]{ALSstring}. Below, we calculate a Weyl function for our example. \\

In view of \eqref{E14+}, we obtain
\begin{align} & \label{E17}
\ze_1(\la)\ze_2(\la)=1, \quad  \ze_1(\la)+ \ze_2(\la)=-2-\frac{1}{\la}.
\end{align}
Simple calculations using \eqref{E15}  show that $|z_2(\la)|>|z_1(\la)|$ for $\la=\ve \I-\frac{1}{4}$, where $\ve$ is small and positive $(\la \in\BC_+)$.
Since $\ze_1\ze_2=1$,  the equality $$\ze_1(\la)=\ov{\ze_2(\la)}$$ is 
valid in  the points $\la\in \BC$, where $|\ze_1(\la)|=1$ or, equivalently, $|\ze_2(\la)|=1$. Hence, in these points $\ze_1+\ze_2\in \BR$. Now, the relation $\ze_1+ \ze_2=-2-\frac{1}{\la}$ implies that
$|\ze_i(\la)|\not=1$ for $\la\in \BC_+$. These considerations show that $|\ze_1(\la)|<1$ for $\la \in \BC_+$. Therefore,
the pair  
\begin{align} & \label{E18}
\clp_1\equiv \E^{-\I c r}  I_p,\quad  \clp_2=\E^{\I c r}\ze_1(\la)\a
\end{align}
is nonsingular, with property-$j$.  From the second equality in \eqref{E14} and relations \eqref{E16} and \eqref{E18}, we derive
\begin{align} & \label{E19}
\clw(r,\la)\begin{bmatrix}\clp_1(r,\la) \\ \clp_2(r,\la)\end{bmatrix}=\E^{-\I \ze_1(\la) r} \begin{bmatrix}I_p \\ \ze_1(\la)\a \end{bmatrix}.
\end{align}
According to \eqref{P7} and \eqref{E19}, we have $ \ze_1(\la)\a \in \cln(r)$ for any $r>0$, that is,
\begin{align} & \label{E20}
\vp(\la)= \ze_1(\la)\a \in \bigcap_{r>0}\cln(r).
\end{align}
Taking into account Definition \ref{DnW}, we see that $ \ze_1(\la)\a$ is a Weyl function of the canonical system \eqref{1.2}, \eqref{1.3}
$($in the case \eqref{E10}, \eqref{E13}$)$ on $[0,\infty)$.
\end{Ee}
A simpler example is considered in the present  paper as Example \ref{Ee6.5}.

Finally, we have the following proposition.
\begin{Pn} \label{PnH} Let generalised canonical system \eqref{1.2} be given on $[0,\infty)$. Then, the sets 
$\cln(r)$ are nested $($i.e., $\cln(r_2)\subseteq \cln(r_1)$ for $0\leq r_1<r_2)$.
Morever, there is a matrix function $\vp(\la)$ belonging to the intersection of all $\cln(r)$, that is,
the set of Weyl functions of the system  \eqref{1.2} on $[0, \, \infty)$ $($the set of $\vp$ satisfying \eqref{P8}$)$ is nonempty.

If \eqref{P8} holds, then the inequality 
\begin{align} & \label{P5}
\int_0^{\infty}\begin{bmatrix}I_{m_1} & \vp(\la)^*\end{bmatrix}W(x,\la)^*H(x)W(x,\la)\begin{bmatrix}I_{m_1} \\ \vp(\la)\end{bmatrix}dx<\infty \quad (\la \in \BC_+)
\end{align}
 is valid.
\end{Pn}
We note that \eqref{P5} may be used as an alternative definition of the Weyl function (see, e.g., \cite{ALSstring}).

{\bf 2.} Similar to  \cite[Appendix A]{ALSstring}, the formulas, statements and proofs  in \cite[Appendix C]{ALSstring} remain
valid after we switch to the case corresponding to a more general form of $j$ (more precisely, when we switch to $m_2\times m$ matrix functions $\b(x)$). 
The corresponding results are given below. 

\begin{Tm}\label{TmSim} Let  the $m_2 \times m$ matrix function $\b(x)$ belong to the class $\clu^{m_2\times m}[0,r]$ $($defined in \eqref{1.3+}$)$
and let \eqref{1.4} hold. Introduce operators $A$ and $K$ acting in $L_2^{m_2}(0,r)$ by the equalities
\begin{align} & \label{P9}
A:=\int_0^x(t-x) \cdot dt, \quad K:=\I \b(x)j\int_0^x\b(t)^*\cdot dt.
\end{align}
Then,  $K$ is linear similar to $A\, :$
\begin{align} & \label{AC3}
K=VAV^{-1}, \quad V=u(x)\big(I_{m_2}+\int_0^x \clv(x,t)\, \cdot \, dt\big),
\end{align}
where $u\in \clu^{m_2\times m_2}[0,r]$,  $u^*=u^{-1}$, and   
\begin{align} & \label{AC4}
\sup\|\clv(x,t)\|<\infty \quad (0\leq t\leq x\leq r).
\end{align}
\end{Tm}
\begin{Rk}\label{RkLT} It is important  for the study of the generalised canonical systems \eqref{1.2}--\eqref{1.4} on the semi-axis $[0, \infty)$
that the matrix function
$\clv(x,t)$ in the domain $0\leq t\leq x\leq \ell$ is uniquely determined by
$\b(x)$ on $[0,\ell]$ $($and does not depend on the choice of $\b(x)$ for $\ell<x<r$ and the choice of $r\geq \ell)$.
\end{Rk}
\section{Fundamental solution}\label{Fund}
\setcounter{equation}{0}
The following considerations are similar to the considerations of \cite[Section~2]{ALSstring} although we consider
a more general case of $j$ and normalised transformation operators $E$ given by \eqref{A36} instead of the
transformation operators~$V$.
Recalling definition \eqref{P9}, it is easy to see that
\begin{align} & \label{F2}
K-K^*=\I \b(x)j\int_0^{ r}\b(t)^* \cdot dt.
\end{align}
If $\b^{\prime\prime}(x)\in L_2^{m_2\times m}(0, \, r)$, we have (according to Proposition \ref{PnM}) $K=EAE^{-1}$, which we substitute
into \eqref{F2}. Multiplying both parts of the derived equality by $E^{-1}$ from the left  and by  $(E^*)^{-1}$ from the right, 
we obtain the operator identity
\begin{align} & \label{F3}
AS-SA^*=\I\Pi j \Pi^*,
\end{align}
where
\begin{align} & \label{F4}
S=E^{-1}(E^*)^{-1}>0, \quad \Pi h=\Pi(x)h, \quad \Pi(x):=\big(E^{-1}\b\big)(x), \\
& \label{F5}
\Pi\in \bB\big(\BC^{m},\, L_2^{m_2}(0, r)\big), \quad  \Pi(x) \in \clu^{m_2\times m}[0,\, r],
\quad h\in \BC^{m}.
\end{align}
Note that $\Pi$ above is the operator of multiplication by the matrix function $\Pi(x)$ and the operator $E^{-1}$ is applied  to $\b$ (in the expression
$E^{-1}\b$) columnwise. The transfer matrix function corresponding  to the so called $S$-node (i.e., to the triple $\{A,S,\Pi\}$ satisfying \eqref{F3})
has the form
\begin{align} & \label{F6}
w_A(\la)=w_A(r,\la)=I_{m}-\I j\Pi^*S^{-1}(A-\la I)^{-1}\Pi,
\end{align}
and was first introduced and studied in \cite{SaL1}. We introduce the projectors $P_{\ell}\in \bB\big(L_2^{m_2}(0,\, r), \, L_2^{m_2}(0,\ell)\big)$:
\begin{align} & \label{F7}
\big(P_{\ell}f\big)(x)=f(x) \quad (0 < x < \ell, \quad \ell \leq r).
\end{align}
Now, we set
\begin{align} & \label{F8}
S_{\ell}:=P_{\ell}SP_{\ell}^*, \quad E_{\ell}:=P_{\ell}EP_{\ell}^*, \quad A_{\ell}:=P_{\ell}AP_{\ell}^*, \quad \Pi_{\ell}:=P_{\ell}\Pi,
\\ & \label{F9}
w_A(\ell,\la)=I_{m}-\I j\Pi_{\ell}^*S_{\ell}^{-1}(A_{\ell}-\la I)^{-1}\Pi_{\ell}.
\end{align}
Since $E$ is a triangular operator, $E^{-1}$ is triangular as well (see, e.g., the formula \eqref{A43} and its proof), and we have $P_{\ell}E^{-1}= P_{\ell}E^{-1}P_{\ell}^*P_{\ell}$.
Hence, taking into account \eqref{F4} and \eqref{F8} we derive
\begin{align} & \label{F10}
P_{\ell}E^{-1}P_{\ell}^*E_{\ell}=P_{\ell}E^{-1}P_{\ell}^*P_{\ell}EP_{\ell}^*=P_{\ell}E^{-1}EP_{\ell}^*=I,
\\ & \label{F11}
S_{\ell}=P_{\ell}E^{-1}(E^*)^{-1}P_{\ell}^*=P_{\ell}E^{-1}P_{\ell}^*P_{\ell}(E^*)^{-1}P_{\ell}^*.
\end{align}
It follows that
\begin{align} & \label{F12}
E_{\ell}^{-1}=P_{\ell}E^{-1}P_{\ell}^*, \quad S_{\ell}=E_{\ell}^{-1}(E_{\ell}^*)^{-1}.
\end{align}
We also have $P_{\ell}A= P_{\ell}AP_{\ell}^*P_{\ell}$. Thus, multiplying both parts of  \eqref{F3} by $P_{\ell}$ from the left and by  $P_{\ell}^*$ from
the right (and using \eqref{F8}, \eqref{F12}, and the last equality in \eqref{F4}) we obtain
\begin{align} & \label{F13}
A_{\ell}S_{\ell}-S_{\ell}A_{\ell}^*=\I\Pi_{\ell} j \Pi_{\ell}^*, \quad \Pi_{\ell}(x)=\big(E_{\ell}^{-1}\b\big)(x) \quad (0<x<\ell).
\end{align}
Clearly $w_A(\ell,\la)$ coincides with $w_A(r, \la)$ when $\ell=r$.
\begin{Rk}\label{wA}
Relations \eqref{F9}, \eqref{F12} and \eqref{F13} show that $S_{\ell}$ and $w_A(\ell,\la)$ may be defined via $E_{\ell}$ $($and $\b(x)$ given on $[0,\ell])$
precisely in the same way as 
$w_A(r,\la)$ is constructed via $E$ $($and $\b(x)$ given on $[0,r])$. Moreover, according to Remark \ref{RkALT}, $E_{\ell}$ may be constructed in the same way as $E$, and so $w_A(\ell,\la)$ does not depend on the choice of $\b(x)$ for $\ell<x<r$ and the choice of $r\geq \ell$. In particular, $w_A(\ell,\la)$ is uniquely defined on the semi-axis
$0<\ell<\infty$ for $\b(x)$  considered on the semi-axis $0\leq x<\infty$.
\end{Rk}
The fundamental solution of the canonical system \eqref{1.2}--\eqref{1.4}
 may be expressed via the  transfer functions $w_A(\ell,\la)$ using continuous
factorisation theorem \cite[p. 40]{SaL2} (see also \cite[Theorem 1.20]{SaSaR} as a more convenient
for our purposes presentation).
\begin{Tm} \label{TmFundSol} Let the Hamiltonian of the generalised canonical system \eqref{1.2} have the form \eqref{1.3}. Assume that $\b(x)$ in \eqref{1.3}
belongs $\clu^{m_2\times m}[0,r]$ and satisfies \eqref{1.4}.
Then, the fundamental solution $W(\ell,\la)$ of the generalised canonical system normalised by  \eqref{P1} admits representation
\begin{align} & \label{F14}
W(\ell,\la)=w_A\Big(\ell,\frac{1}{\la}\Big) 
\end{align}
for $ 0<\ell\leq r$. If theorem's conditions hold for each $0<r<\infty$, then \eqref{F14} is valid for each $\ell$ on the semi-axis $(0,\infty)$.
\end{Tm}
The proof of Theorem \ref{TmFundSol} coincides with the proof of  Theorem 2.2 in \cite{ALSstring}.
\begin{Rk}\label{RkPhi2} Using \eqref{F4}, \eqref{F5} and the last equality in  \eqref{A42}, we  partition $\Pi$ and $\Pi(x)$ into two blocks and derive:
\begin{align} & \label{F15-}
\Pi=\begin{bmatrix}\Phi_1 & \Phi_2\end{bmatrix}, \quad \Phi_k h_k=\Phi_k(x)h_k \,\, (k=1,2);
\\ &  \label{F15}
\Pi(x)=\begin{bmatrix}\Phi_1(x) & \Phi_2(x)\end{bmatrix}, \quad \Phi_2(x)\equiv I_{m_2}.
\end{align}
\end{Rk}
After partitioning $\b$ into blocks $\b_k\in \clu^{m_2\times m_k}[0,r]$ $(k=1,2)$, relations \eqref{F4} and \eqref{F15} yield
\begin{equation}  \label{F16}
  \Phi_1(x)=\big(E^{-1}\b_1\big)(x), \quad \b(x)=:\begin{bmatrix}\b_1(x) & \b_2(x)\end{bmatrix}.
\end{equation}

\section{Structured operators $S$} \label{Struct}
The study of the structured operators in direct and inverse spectral problems goes back to  the classical  works \cite{Krein0, Krein} by M.G. Krein.
See further developments, discussions and references in \cite{GeS,  ALS02, SaSaR, SaL2-, LA94, SaL2}.
The case of $S$ satisfying  operator identity (or equation) \eqref{F3}, where $A$ is given in \eqref{P9}, that is, 
\begin{align} & \label{S1}
A=\cla^2, \quad \cla=\I \int_0^x\cdot \, dt \in \bB\big(L_2^{m_2}(0,r)\big),
\end{align}
was studied only in \cite{DuK}. However, related Theorems 1.1 and 1.2 in \cite{DuK} contain mistakes (since the functions in (1.4) and (1.18) in \cite{DuK}
are integrated over domains, where they are not defined).  Hence, we will  consider  \eqref{F3} here. In particular, we will
study $S$ satisfying an important operator identity
\begin{align} & \label{S2}
\cla S+S\cla^*=R, \quad R:=\int_0^r\clr(x,t) \cdot dt \in \bB\big(L_2^{m_2}(0,r)\big).
\end{align}
Contrary to the identity $\cla S-S\cla^*=\I R$, there is no literature on \eqref{S2} as far as we know.

Since $A=\cla^2$, the operator $AS-SA^*$ may be rewritten in the form
\begin{align} & \label{S3}
AS-SA^*=\cla(\cla S+S\cla^*)-(\cla S+S\cla^*)\cla^*,
\end{align}
that is,  \eqref{F3} may be rewritten as
\begin{align} & \label{S3+}
\cla Z-Z\cla^*=\I \Pi j \Pi^*, \quad \cla S+S\cla^*=Z.
\end{align}
Thus, operators $S$ satisfying \eqref{F3} may be found in two steps, which we discuss below.
First, we study $Z$ such that
\begin{align} & \label{S4}
\cla Z-Z\cla^*=\I \Pi j \Pi^*.
\end{align}
The operators $Z$ satisfying  \eqref{S4} as well  as more general identity 
\begin{align} & \label{S4+}
\cla Z-Z\cla^*=\I R, \quad R:=\int_0^r\clr(x,t) \cdot dt \in \bB\big(L_2^{m_2}(0,r)\big)
\end{align}
were thoroughly studied in \cite{SaL1+, SaL15} (see also \cite{KKL}).  According to the resulting
Theorem D.1 \cite{SaSaR}, there is a unique solution $Z \in \bB\big(L_2^{m_2}(0,r)\big) $ of the identity \eqref{S4} where $\Pi$ has the form \eqref{F15-}, \eqref{F15}
(such that $\Phi_1(x)$ has a bounded derivative):
\begin{align}&
\label{S5}
\big(Zf\big)(x)=\big(\Phi_1(0)\Phi_1(0)^*-I_{m_2}\big)f(x)+\int_0^r \clz(x,t)f(t)dt, \\
&\label{S6}
\clz(x,t):=\int_0^{\min(x,t)}\Phi_1^{\prime}(x-\zeta)\Phi_1^{\prime}(t-\zeta)^*d\zeta
+
\begin{cases}\Phi_1^{\prime}(x-t)\Phi_1(0)^*, \quad x>t;
\\
\Phi_1(0)\Phi_1^{\prime}(t-x)^*, \quad t>x.
\end{cases}
\end{align}
\begin{Rk}\label{RkGen}
According to \cite[Corollary 1.1.7]{SaL15} the bounded operator $Z$, which satisfies \eqref{S4+} $($that is,  a more general identity
than \eqref{S4}$)$, admits
representation
\begin{align} & \label{Z1}
Z=\frac{d}{dx}\int_0^r\left(\frac{\p}{\p t}\Psi(x,t)\right) \, \cdot \, dt,
\\ & \label{Z2}
\Psi(x,t) =\frac{1}{2}\int_{|x-t|}^{x+t}\clr\left(\frac{s+x-t}{2},\frac{s-x+t}{2}\right)ds.
\end{align}
\end{Rk}

In our case,  \eqref{F16} and Proposition \ref{PnM} imply that $\Phi_1(x)\in \clu^{m_2\times m_1}[0,r]$, and so the relations
\eqref{S5} and \eqref{S6} hold. Moreover, taking into account \eqref{F16}, \eqref{A8!} and \eqref{A36+} we obtain
\begin{align} & \label{S7}
\Phi_1(0)=\b_2(0)^{-1}\b_1(0),
\end{align}
where $\b_2(0)$ is invertible.
In view of the first equality in \eqref{1.4} at $x=0$, formula \eqref{S7} yields 
\begin{align} & \label{S8}
\Phi_1(0)\Phi_1(0)^*=I_{m_2}.
\end{align}
Therefore, we rewrite \eqref{S5} in a simpler form
\begin{align} & \label{S9}
Z=\int_0^r\clz(x,t)\, \cdot \, dt \in \bB\big(L_2^{m_2}(0,r)\big).
\end{align}
Next, we consider a similar to \eqref{S2} identity
\begin{align} & \label{S10}
 \cla  S+S \cla^*=Z.
\end{align}
In order to study the operator identity \eqref{S10}, we start with the closely related identity
\begin{align} & \label{S11}
T \cla +\cla^*T=\wt Z, \quad \wt Z=\int_0^r \wt \clz(x,t)\, \cdot \, dt \in \bB\big(L_2^{m_2}(0,r)\big),
\end{align}
and modify the proof of \cite[Theorem 1.1.6]{SaL15} (i.e., of Theorem 1.3 in \cite[Chapter 1]{SaL1+}).
\begin{Tm}\label{TmT} Let operator $T$ belong to the class $ \bB\big(L_2^{m_2}(0,r)\big)$ of bounded operators and satisfy operator identity \eqref{S11}.
Then, $T$ admits representation
\begin{align} & \label{S12}
T=\frac{d}{dx}\int_0^r\frac{\p \wt\Up}{\p t}(x,t) \, \cdot \, dt,
\end{align}
where $\wt \Up(x,t)$  has the form
\begin{align} & \label{S13}
\wt\Up(x,t)=\wt \up(x+t)+\frac{\I}{2}\int_{-\min\{x+t,\, 2r-x-t\}}^{x-t}\wt \clz\left(\frac{x+t+s}{2}, \frac{x+t-s}{2}\right)ds,
\\ & \label{S13+}
\frac{\p}{\p t}\wt\Up(x,t) \in L_2^{m_2\times m_2}(0,r)  \quad  {\mathrm{for}} \quad x\in [0,r], 
\\ & \label{S14}
\wt \up(x)\in L_2^{m_2\times m_2}(0,r), \quad \wt \up(x)=0 \quad {\mathrm{for}} \quad r\leq x\leq2r.
\end{align}
Moreover, in this case the kernel $\wt \clz(x,t)$ satisfies the condition
\begin{align} & \label{S15}
\int_{\xi-r}^r\wt \clz(\xi-t,t) dt=0 \quad  {\mathrm{for}} \quad  r\leq \xi \leq 2r.
\end{align}
\end{Tm}
\begin{proof}.  Representation \eqref{S12} of the bounded operators $T$ follows from \cite[Theorem 1.1.1]{SaL15}.
Clearly, we may assume additionally that
\begin{align} & \label{S15!}
\wt \Up(r,t)=0.
\end{align}
It remains to prove relations \eqref{S13}, \eqref{S14} and \eqref{S15}. For this purpose, one can use the first part
of the proof of \cite[Theorem 1.1.6]{SaL15}, where the signs of the second terms on the left-hand sides of  (1.1.25), (1.1.28)
and (1.1.29) in
\cite{SaL15} should be (evidently) changed.  We also change some notations in \cite{SaL15}: $Q_2$ into $\clz_2$, $\I Q$ into  $\wt \clz$ and $\om$ into $r$.
Then,  in view of  \cite[(1.1.22)]{SaL15} and of the second formula in \cite[(1.1.32]{SaL15}
we have
\begin{align}
& \label{S17}
 F_1(x,t):=-\int_x^r\int_t^r \frac{\p \wt \Up}{\p \zeta}(y,\zeta) \, d\zeta \, dy, \quad Z_2(x,t):=-\I \int_x^r \wt \clz(y,t)dy.
\end{align}
Taking into account the abovementioned changes, instead of  \cite[(1.1.30)]{SaL15}
we  obtain the following partial differential equation for $F_1$:
\begin{align} & \label{S16}
\frac{\p F_1 }{\p t}(x,t)-\frac{\p F_1}{\p x}(x,t)=Z_2(x,t).
\end{align}
Equation  \eqref{S16} may be rewritten in the form
\begin{equation}  \label{S18}
-2\frac{\p F_2 }{\p \eta}(\xi,\eta)=\clz_2\left(\frac{\xi+\eta}{2}, \frac{\xi-\eta}{2}\right), \quad F_2(\xi, \eta):=F_1\left(\frac{\xi+\eta}{2}, \frac{\xi-\eta}{2}\right).
\end{equation}
According to \eqref{S17}, we have
\begin{align} & \label{S19}
F_1(r,t)=F_1(x,r)=0,
\end{align}
which yields
\begin{align} & \label{S20}
F_2\big(\xi,\pm(2r-\xi)\big)=0 \quad  {\mathrm{for}} \quad  r\leq \xi \leq 2r.
\end{align}
Introduce the boundary value matrix function
\begin{equation}  \label{S21}
\vt(\xi)=\begin{cases} F_2(\xi,-\xi)=F_1(0,\xi) \quad  {\mathrm{for}} \quad  0\leq \xi < r, \\ 0 \quad  {\mathrm{for}} \quad  r\leq \xi \leq 2r.
\end{cases}
\end{equation}
Using \eqref{S18}, \eqref{S20} and \eqref{S21}, we derive
\begin{align} & \label{S22}
F_2(\xi,\eta)=\vt(\xi)-\frac{1}{2}\int_{-\xi}^{\eta}\clz_2\left(\frac{\xi+y}{2}, \frac{\xi-y}{2}\right)dy \quad  {\mathrm{for}} \quad  0\leq \xi < r,
\\  & \label{S23}
F_2(\xi,\eta)=\vt(\xi)-\frac{1}{2}\int_{\xi-2r}^{\eta}\clz_2\left(\frac{\xi+y}{2}, \frac{\xi-y}{2}\right)dy \quad  {\mathrm{for}} \quad  r\leq \xi \leq 2r.
\end{align}
Taking into account the second equality in \eqref{S18} and setting $\xi=x+t, \,\, \eta=x-t$, we rewrite  \eqref{S22} and \eqref{S23} as
\begin{align} & \label{S24}
F_1(x,t)=\vt(x+t)-\frac{1}{2}\int_{-\min\{x+t, \, 2r-x-t\}}^{x-t}\clz_2\left(\frac{x+t+y}{2}, \frac{x+t-y}{2}\right)dy.
\end{align}
Setting $\zeta=(x+t-y)/2$ and using the second equality in \eqref{S17}, we rewrite \eqref{S24} in the form
\begin{align} & \label{S24'}
F_1(x,t)=\vt(x+t)+\I \int^{\min\{x+t, \, r\}}_{t}\int_{x+t-\zeta}^r\wt \clz\left({s}, \zeta\right)ds\,d\zeta.
\end{align}
Hence, we have
\begin{align} & \label{S25}
\frac{\p}{\p x}F_1(x,t)=-\wt \up(x+t)-\I \int^{\min\{x+t, \, r\}}_{t}\wt \clz\left({x+t-\zeta}, \zeta\right)d\zeta,
\\ & \label{S25'}
\wt \up(x)=-\vt^{\prime}(x)-\I\int_0^r\wt \clz\left({s}, x\right)ds \,\, {\mathrm{for}} \,\, x< r, \quad \wt\up(x)=0 \,\, {\mathrm{for}} \,\, r\leq x\leq 2r.
\end{align}
It follows from \eqref{S15!}, \eqref{S17} and \eqref{S25} that
\begin{align} \nn
\wt \Up(x,t) & =-\frac{\p F_1}{\p x}(x,t)=\wt \up(x+t)+\I \int^{\min\{x+t, \, r\}}_{t}\wt \clz\left({x+t-\zeta}, \zeta\right)d\zeta
\\ & \label{S26}
=\wt \up(x+t)+\frac{\I}{2} \int_{-\min\{x+t, \, 2r-x-t\}}^{x-t}\wt \clz\left(\frac{x+t+s}{2}, \frac{x+t-s}{2}\right)ds.
\end{align}
According to \eqref{S21}, \eqref{S25'} and \eqref{S26}, relations \eqref{S13} and \eqref{S14} are valid.

Finally, in order to prove \eqref{S15} we note that the first equality in \eqref{S18} and relations \eqref{S20} imply that
\begin{align} & \label{S27}
\int_{\xi-2r}^{2r-\xi}\int_{\frac{\xi+\eta}{2}}^r\wt \clz\left(x, \frac{\xi-\eta}{2}\right)dx\,d\eta=0.
\end{align}
Setting $t=\frac{\xi-\eta}{2}$ in \eqref{S27}, we rewrite \eqref{S27} in the form 
\begin{align} & \label{S27+}
\int_{\xi-r}^r\int_{\xi-t}^r \wt \clz(x,t)\, dx\, dt=0 \quad  {\mathrm{for}} \quad  r\leq \xi \leq 2r.
\end{align}
Since the right-hand side of \eqref{S27+} equals zero at $\xi=r$, by differentiating \eqref{S27+} with respect to $\xi$
we obtain an equivalent equality, that is, \eqref{S15} follows.
\end{proof}
Introduce operators $U$ such that
\begin{align} & \label{S28}
\big(Uf\big)(x)=\ov{f(r-x)} \quad \big(f\in L_2^{m_2}(0,r)\big).
\end{align}
It is immediate that
\begin{align} & \label{S29}
U\cla U=\cla^*, \quad U^2= I.
\end{align}
Hence, the identity \eqref{S10} is equivalent to the identity
\begin{align} & \label{S30}
USU\cla+\cla^*USU=U\clz U
\end{align}
and Theorem \ref{TmT} yields the following corollary.
\begin{Cy} \label{CyS}
Let operator $S$ belong to the class $ \bB\big(L_2^{m_2}(0,r)\big)$ of bounded operators and satisfy operator identity \eqref{S10}
where $Z$ has the form \eqref{S9}.
Then, $S$ admits representation
\begin{align} & \label{S31}
S=\frac{d}{dx}\int_0^r\left(\frac{\p}{\p t}\Up(x,t) \right)\, \cdot \, dt,
\end{align}
where $ \Up(x,t)$  has the form
\begin{align} & \label{S32}
\Up(x,t)= \up(x+t)-\frac{\I}{2}\int^{\min\{x+t,\, 2r-x-t\}}_{x-t} \clz\left(\frac{x+t+s}{2}, \frac{x+t-s}{2}\right)ds,
\\ & \label{S32+}
\frac{\p}{\p t}\Up(x,t) \in L_2^{m_2\times m_2}(0,r)  \quad  {\mathrm{for}} \quad x\in [0,r], 
\\ & \label{S33}
  \up(x)=0 \quad {\mathrm{for}} \quad 0\leq x\leq r, \quad \up(x)\in L_2^{m_2\times m_2}(r,2r).
\end{align}
Moreover, in this case the kernel $\clz(x,t)$ satisfies the condition
\begin{align} & \label{S34}
\int_{0}^{\xi}  \clz(\xi-t,t) dt=0 \quad  {\mathrm{for}} \quad  0\leq \xi \leq r.
\end{align}
\end{Cy}
\begin{proof}. Setting 
\begin{align} & \label{S35}
S=UTU, \quad Z=U\wt Z U,
\end{align}
where $T$, $\wt Z$ satisfy \eqref{S11}, and using \eqref{S9} and \eqref{S12} we derive
\eqref{S31}, where
\begin{align} & \label{S36}
\Up(x,t)=\ov{\wt \Up(r-x,r-t)}, \quad \ov{\wt \clz(x,t)}=\clz(r-x,r-t).
\end{align}
According to \eqref{S13} and \eqref{S36} we have
\begin{align} \nn
\Up(x,t)=&\ov{\wt \up(2r-x-t)}
\\ & \nn
-\frac{\I}{2}\int^{t-x}_{-\min\{2r-x-t,\, x+t\}}\ov{\wt \clz\left(\frac{2r-x-t+s}{2}, \frac{2r-x-t-s}{2}\right)}ds
\\ \nn
=&\ov{\wt \up(2r-x-t)}
\\ & \label{S37}
-\frac{\I}{2}\int^{t-x}_{-\min\{2r-x-t,\, x+t\}} \clz\left(\frac{x+t-s}{2}, \frac{x+t+s}{2}\right)ds.
\end{align}
Thus, we obtain \eqref{S32} and \eqref{S33} after we switch from $s$ to $-s$ on the right-hand side
of \eqref{S37}. Finally, in view of the second equality in \eqref{S36}, condition \eqref{S15} may be
rewritten in the form 
\begin{align} & \label{S37+}
\int^{r}_{\xi-r}  \clz(r-\xi+t,r-t) dt=\int^{2r-\xi}_{0}  \clz(2r-\xi-s,s) ds=0, \quad    r\leq \xi \leq 2r.
\end{align}
Substituting in \eqref{S37+} $\xi$ instead of $2r-\xi$, we derive \eqref{S34}.
\end{proof}
\begin{Rk}\label{OpId} Since operators $S_{\ell}$ in \eqref{F9} $($and, correspondingly, in the representation \eqref{F14}
of the fundamental solution$)$ satisfy operator identities  \eqref{F13}, they have the form
\eqref{S31}--\eqref{S33}, where $r$ is substituted by $\ell$ and $\clz$ is given by \eqref{S6}.
\end{Rk}
In order to consider our conditions on $\Phi_1$, it is convenient to write down the right-hand side
of \eqref{S3} in another form:
\begin{align} & \label{S38}
AS-SA^*=\cla(\cla S-S\cla^*)+(\cla S-S\cla^*)\cla^*.
\end{align}
Then, taking into account \eqref{F3} and \eqref{F15}, we rewrite condition \eqref{S34} in the form
\begin{align} & \label{S39}
\int_{0}^{\xi}  \Phi_1(\xi-t)  \Phi_1(t)^* dt=\xi I_{m_2}.
\end{align}
\begin{Cy} \label{CyCond} If the conditions of Theorem \ref{TmFundSol} hold on $[0,r]$ $\big($on $[0,\infty)\big)$,  then \eqref{S39}
is valid for all $\xi \in [0,r]$ $\big(\xi \in [0,\infty)\big)$ as well. 
\end{Cy}
Finally, let us consider some examples.
\begin{Ee}\label{Ee1} Let $Z$ be a bounded operator, which satisfies \eqref{S4+}, where
\begin{align} & \label{Z3}
\clr(x,t)=\clr_0(x+t), \quad \clr_0^{\prime}(s)\in L_2^{m_2\times m_2}(0,2r).
\end{align}
Then, \eqref{Z2} yields
\begin{align} & \label{Z4}
\Psi(x,t)=\frac{1}{2}\int_{|x-t|}^{x+t}\clr_0(s)ds,
\end{align}
and explicit representation of $Z$ follows from \eqref{Z1}, \eqref{Z4}$:$ 
\begin{align} & \label{Z5}
\big(Zf\big)(x) =\clr_0(0)f(x)+\frac{1}{2}\int_0^r\Big(\clr_0^{\prime}(x+t)+\clr_0^{\prime}(|x-t|)\Big)f(t) dt.
\end{align}
\end{Ee}

It is also easily checked that each $Z$ of the form \eqref{Z5} (where $\clr_0$ is differentiable and
$\clr_0^{\prime}(s)\in L_2^{m_2\times m_2}(0,2r)$) satisfies the  identity \eqref{S4+} with
\begin{align} & \label{Z6}
R= \int_0^r\clr_0(x+t) \, \cdot \, dt.
\end{align}
Indeed, we set
\begin{align} & \nn
Z=Z_1+Z_2+Z_3, \quad Z_1f=\clr_0(0)f, \quad Z_2f=\frac{1}{2}\int_0^r\clr_0^{\prime}(x+t)f(t) dt,
\\ & \nn
 Z_3f=\frac{1}{2}\int_0^r\clr_0^{\prime}(|x-t|)f(t) dt.
\end{align}
Clearly, 
\begin{align} & \label{Z7}
(\cla Z_1-Z_1\cla^*)f=\I\int_0^r f(t)dt.
\end{align}
Simple calculations similar to \eqref{Z21} and \eqref{Z22} below show that
\begin{align} & \label{Z8}
(\cla Z_2-Z_2\cla^*)f=\frac{\I}{2}\int_0^r\big(2\clr_0(x+t)-\clr_0(t)-\clr_0(x)\big)f(t)dt.
\end{align}
Finally, changing the order of integration, we obtain
\begin{align}\nn
\cla Z_3f=&\frac{\I}{2}\left(\int_0^x\big(\clr_0(x-t)+\clr_0(t)-2\clr_0(0)\big)f(t)dt \right.
\\  & \label{Z9}
\left. +\int_x^r\big(\clr_0(t)-\clr_0(t-x)\big)f(t)dt\right),
\\ \nn
Z_3\cla^*=&-\frac{\I}{2}\left(\int_0^x\big(\clr_0(x)-\clr_0(x-t)\big)f(t)dt\right.
\\  & \label{Z10}
\left. +\int_x^r\big(\clr_0(x)+\clr_0(t-x)-2\clr_0(0)\big)\right).
\end{align}
Equalities \eqref{Z7}--\eqref{Z10} imply the identity \eqref{S4+} where $R$ has the form \eqref{Z6}.
\begin{Ee}\label{Ee2} The operator
\begin{align} & \label{S40}
S_0=\int_0^rv_0(x+t)\,\cdot \, dt  \quad \big(v_0\in L_2^{m_2\times m_2}(0,2r)\big)
\end{align}
satisfies the operator identity
\begin{align} & \label{S41}
\cla S_0+S_0\cla^*=\I\int_0^r\int_t^x v_0(s)ds \, \cdot \, dt.
\end{align}
\end{Ee}
Indeed, we have
\begin{align} & \label{Z21}
\big(\cla S_0f\big)(x)=\I\int_0^x\int_0^rv_0(t+s)f(s)dsdt=\I\int_0^r\int_0^xv_0(t+s)dtf(s)ds, \\
& \label{Z22}
\big(S_0\cla^*f\big)(x)=-\I \int_0^rv_0(x+t)\int_t^rf(s)dsdt=-\I\int_0^r\int_0^sv_0(x+t)dtf(s)ds,
\end{align}
and \eqref{S41} follows.

According to \eqref{S32}, \eqref{S33} and \eqref{S41}, $\Up(x,t)$ is given in this case by the equalities
\begin{align} & \label{S42}
\Up(x,t)=\frac{1}{2}\int_{x-t}^{x+t}\int_{(x+t-s)/2}^{(x+t+s)/2}v_0(\xi)d\xi d s \quad (x+t\leq r),
\\ & \label{S43}
\Up(x,t)=\up(x+t)+\frac{1}{2}\int_{x-t}^{2r-x-t}\int_{(x+t-s)/2}^{(x+t+s)/2}v_0(\xi)d\xi d s \quad (x+t\geq r).
\end{align}
Hence, we derive
\begin{align} & \label{S44}
\frac{\p}{\p t}\Up(x,t)=\int^{x+t}_t v_0(\xi)d\xi \quad (x+t\leq r), \\
 & \label{S45}
 \frac{\p}{\p t}\Up(x,t)=\up^{\prime}(x+t)-\int_{x+t-r}^tv_0(\xi)d\xi \quad (x+t\geq r).
\end{align}
Comparing \eqref{S31} and \eqref{S40} and taking into account \eqref{S44} and \eqref{S45},
we see that
\begin{align} & \label{S46}
\up^{\prime\prime}(\xi)=v_0(\xi)-v_0(\xi-r) \quad {\mathrm{for}} \quad \xi>r
\end{align}
in our example. The condition \eqref{S34} takes the form
\begin{align} & \label{S47}
\I\int_0^\xi\int_t^{\xi-t}v_0(s)dsdt=0 \quad (0\leq \xi \leq r),
\end{align}
which is equivalent to
\begin{align} & \label{S48}
\left(\int_0^\xi\int_t^{\xi-t}v_0(s)dsdt\right)^{\prime}=0 \quad (0\leq \xi \leq r).
\end{align}
Clearly, \eqref{S48} is always valid.

Using our results for Examples \ref{Ee1} and \eqref{Ee2} we obtain one more example. 
\begin{Ee}\label{Ee3} Let $S=Z$ be given by \eqref{Z5} and let $A= \cla^2$
$($i.e., let $A$ be given by  \eqref{P9}$)$. Then, $S$ satisfies the operator
identity 
\begin{align} & \label{S49}
AS-SA^*=-\int_0^r\int_t^x\clr_0(s)ds \,\cdot \, dt.
\end{align}
Moreover, the rank of the right-hand side of \eqref{S49} is no more than $2m_2$.
\end{Ee}

Indeed, for $S_0$ given \eqref{S40}, relations \eqref{S4+} and \eqref{Z6}  yield
\begin{align} & \label{S50}
\cla S-S\cla^*=S_0, \quad {\mathrm{where}}  \quad v_0(s):=\I \clr_0(s).
\end{align}
Now, using equalities \eqref{S41} and \eqref{S50} we obtain
\begin{align} & \label{S51}
\cla(\cla S-S\cla^*)+(\cla S-S\cla^*)\cla^*=\cla S_0+S_0\cla^*=-\int_0^r\int_t^x\clr_0(s)ds \,\cdot \, dt.
\end{align}
 Since the left-hand sides of \eqref{S49} and \eqref{S51} coincide, the operator  identity \eqref{S49} follows.
 The estimate of the rank of the right-hand side of \eqref{S49} follows from the representation
 $$\int_t^x\clr_0(s)ds=\int_0^x\clr_0(s)ds-\int_0^t\clr_0(s)ds,$$
where the first and second terms on the right-hand side are  $m_2\times m_2$ matrix functions, the first one depending
only on $x$ and the second one depending on $t$.
 
Operators $S$ of the form \eqref{Z5}, some similar ones, and spectral theory of the corresponding string equations
have been studied, for instance, in \cite{Krein0, Langer} (in the scalar case), in \cite[pp. 53-56]{SaL2-} as well as in \cite[\S 7]{ALS88}.

In the spirit of the continuous factorisation theorem \cite{SaSaR, SaL2-, SaL2}, formula \eqref{F4} and Remark \ref{wA} yield
$\Pi_r^*S_r^{-1}\Pi_r=\int_0^r\b(x)^*\b(x)dx$, that is,
\begin{align} & \label{S52}
\frac{d}{dr}\Big(\Pi_r^*S_r^{-1}\Pi_r\Big)=\b(r)^*\b(r)=H(r).
\end{align}
The recovery of $\Phi_1$ and operators $S_r$, and so (in view of \eqref{S52}) the solution of the inverse problem to
recover a generalised canonical system, is the subject of our next paper.
\section{High energy asymptotics \\ of the Weyl functions} \label{Direct}
\setcounter{equation}{0}
Relations \eqref{F3}, \eqref{F6} and \eqref{F14} yield (see \cite[(1.88)]{SaSaR}):
\begin{align} & \label{M1}
W(r,\la)^*jW(r,\la)=j+\I(\la -\ov{\la})\Pi^*(I-\ov{\la}A^*)^{-1}S^{-1}(I-{\la}A)^{-1}\Pi,
\end{align}
where $A=A_r, \, S=S_r,\, \Pi=\Pi_r$. Hence, in view of \eqref{P7+} we obtain
\begin{align} & \label{M2}
I_{m_1}\geq \I(\ov{\la}-\la)\begin{bmatrix}I_{m_1} & \phi(\la)^*\end{bmatrix}\Pi^*(I-\ov{\la}A^*)^{-1}S^{-1}(I-{\la}A)^{-1}\Pi\begin{bmatrix}I_{m_1} \\ \phi(\la)\end{bmatrix}
\end{align}
for $\la \in \BC_+$ and $\phi(\la)\in \cln(r)$. Since the operator $S$ is strictly positive, \eqref{M2} implies that
\begin{equation}  \label{M3}
\left\|(I-\la A)^{-1}\Pi\begin{bmatrix}I_{m_1} \\ \phi(\la)\end{bmatrix}\right\| \leq \frac{ C}{\sqrt{\Im(\la)}} \quad {\mathrm{for}} \,\, \la \in \BC_+ \,\,
{\mathrm{and \,\, some}} \quad C=C(r)>0.
\end{equation}
It is easily checked directly and follows from the formula for the resolvent 
$(I+z\cla)^{-1}=I-\I z\int_0^x\E^{\I z(t-x)}\, \cdot\, dt$
 (see, e.g., \cite[(1.157)]{SaSaR}) that
\begin{align} & \label{M4}
(I-z^2A)^{-1}=I+\frac{\I z}{2}\int_0^x\Big(\E^{\I z(x-t)}-\E^{\I z(t-x)}\Big) \,\cdot \, dt.
\end{align}
In particular, for the operators $\Phi_k\in  \bB\big(\BC^{m_k},\, L_2^{m_2}(0, r)\big)$ $(k=1,2)$, which are
defined in \eqref{F15-} and \eqref{F15}, we have
\begin{align} & \label{M5}
(I-z^2A)^{-1}\Phi_2=\frac{1}{2}\Big(\E^{\I zx}+\E^{-\I z x}\Big)\Phi_2, \\
& \label{M6}
 \Phi_2^*(I-z^2A)^{-1}\Phi_2=\frac{1}{2\I z}\Big(\E^{\I zr}-\E^{-\I z r}\Big)I_{m_2}, \\
&  \label{M7} 
 \Phi_2^*(I-z^2A)^{-1}\Phi_1=\frac{1}{2}\int_0^r\Big(\E^{\I z(r-t)}+\E^{\I z(t-r)}\Big)\Phi_1(t)dt.
\end{align}
It is immediate from \eqref{M6} that
\begin{align} & \label{M8}
\big( \Phi_2^*(I-z^2A)^{-1}\Phi_2\big)^{-1}= -2\I z \E^{\I z r}\big(1+O(\E^{2 \I z r})\big)I_{m_2}\quad {\mathrm{for}} \quad  \Im(z)\to \infty .
\end{align}
In turn, relations \eqref{M3}, \eqref{M7} and \eqref{M8} yield
\begin{equation}  \label{M9}
\phi(z^2)=\left(o(1)+\frac{1}{\sqrt{\Im(z^2)}}\right)|z| O(\E^{ \I z r})+\I z\int_0^r\E^{ \I z t}\Phi_1(t)dt \quad \big( \Im(z)\to \infty\big) .
\end{equation}
Note that $O(\E^{ \I z r})$ in \eqref{M9} characterises the growth of the corresponding matrix norm and $z$ is assumed
to be situated in the first quadrant ($\Re(z)>0$, $\Im z >0$), that is, $\Im(z)>0$ and $\Im(z^2)>0$.
\begin{Tm}\label{TmHEA}
Let generalised  canonical system \eqref{1.2}, \eqref{1.3} be given, such that $\b(x)\in \clu^{m_2\times m}[0,r]$ and \eqref{1.4} holds.
Then, the Weyl functions $\phi~\in~\cln(r)$ admit representation \eqref{M9}, where $\Phi_1(x)$ is given by  \eqref{F16}. The corresponding
matrix function $\Phi_1(x)$
is two times differentiable and belongs to $\clu^{m_2\times m_1}[0,r]$. Moreover, $\Phi_1(x)$ satisfies the equality
\begin{align} & \label{S39'}
\int_{0}^{\xi}  \Phi_1(\xi-t)  \Phi_1(t)^* dt=\xi I_{m_2}.
\end{align}
The  matrix function $\Phi_1(x)$  is uniquely
determined on $[0,r]$ by  the representation \eqref{M9}. 
\end{Tm}
\begin{proof}. Formula \eqref{M9} was proved above.  The equality \eqref{F16} and Proposition \ref{PnM} show that 
$\Phi_1(x)\in \clu^{m_2\times m_1}[0,r]$. Condition \eqref{S39'} follows from Corollary \ref{CyCond}.

Finally, suppose that there is a matrix function $\br \Phi(x)\in L_2^{m_2\times m_1}(0,r)$ such that
\begin{equation}  \label{M10}
\phi(z^2)=\left(o(1)+\frac{1}{\sqrt{\Im(z^2)}}\right)|z| O(\E^{ \I z r})+\I z\int_0^r\E^{ \I z t}\br \Phi(t)dt.
\end{equation}
The equalities \eqref{M9} and \eqref{M10}
yield
\begin{align} & \label{M11}
\om(z):=\int_0^r\E^{ \I z (t-r)}\big(\Phi_1(t)-\br \Phi(t)\big)dt=\left(o(1)+\frac{1}{\sqrt{\Im(z^2)}}\right)O(1).
\end{align}
Clearly, $\| \om(z)\|$ is bounded in the domain $\{z: \, \Im(z)\leq \ve\}$ $(\ve>0)$. According to \eqref{M11},
$\| \om(z)\|$ is bounded in the domain $\{z: \, \Re(z)\geq \ve, \,\, \Im(z)\geq \ve\}$.  Finally, by virtue of the
Phragmen--Lindel\"of theorem (see, e.g., \cite[Corollary E.7]{SaSaR} after using the change of variables $z=(\exp\{3\I \pi/4\})\wt z+\ve(1+\I)$
in order to come to the standard angle considered in the theorem),
$\|\om(z)\|$ is bounded in the domain $\{z: \, \Re(z)< \ve, \,\, \Im(z) > \ve\}$. Thus, $\|\om(z)\|$ is bounded in $\BC$
and $\om(z)$ tends also to zero on some rays in $\BC$. Hence, we obtain $\om(z)\equiv 0$, and the identity $\br \Phi(x) \equiv \Phi_1(x)$
follows.
\end{proof}
\begin{Rk}\label{RkL}
Since $\Phi_1(x)$ is two times differentiable, the integral on the right-hand side of \eqref{M9} can be integrated by parts, which
will produce a slightly more precise asymptotics.
\end{Rk}
\section{A special case $m_1=m_2=p$: Weyl functions from Herglotz class} \label{WH}
\setcounter{equation}{0}
{\bf 1.} In the important case $m_1=m_2=p$, together with system \eqref{1.2} studied above, we also consider the standard form \eqref{1.1} of the canonical system:
\begin{align} &       \label{B1}
\wh w^{\prime}(x,\la)=\I \la J \wh H(x)\wh w(x,\la),
 \end{align} 
 and assume that
 \begin{align}&\label{B2}
\wh H(x)=\wh \b(x)^*\wh \b(x), \quad \wh \b(x)J\wh \b(x)^*\equiv 0, \quad  \wh \b^{\prime}(x)J\wh \b(x)^* \equiv \I I_{p},
 \end{align}
where $J$ is given in \eqref{1.1} and the ``widehat" is used (in the notations of this section) in order to show that the notations
correspond to the system \eqref{1.1} (or, more precisely, \eqref{B1}) instead of the system \eqref{1.2}.
We require that $\wh \b(x)\in \clu^{p\times 2p}[0,r]$. According to \cite[(1.3)]{ALSstring}, we have
\begin{align}& \label{B3}
J=\Theta j \Theta^*, \quad \Theta:=\frac{1}{\sqrt{2}}\begin{bmatrix} I_p & -I_p \\ I_p & I_p\end{bmatrix}  \quad (\Theta \Theta^*=I_{2p}). \end{align}
Thus, one may assume the following simple correspondence between the systems \eqref{1.2}, \eqref{1.3}, \eqref{1.4} (where $m_1=m_2=p$) and 
systems \eqref{B1}, \eqref{B2}:
\begin{align} &       \label{B4}
\wh\b(x)=\b(x)\Th^*, \quad  \wh H(x)=\Th H(x)\Th^*, \quad \wh W(x,\la)=\Th W(x,\la)\Th^*.
 \end{align} 
Here, $\wh W(x,\la)$ is the normalised fundamental solution of \eqref{B1}.  
Similar to $\cln(r)$, the set $\wh \cln(r)$ is introduced as the set of linear-fractional transformations $\wh \phi$:
\begin{align} \nn
\wh \phi(r,\la)=&\I \big(\wh \clw_{21}(r,\la)\wh \clp_1(\la)+\wh \clw_{22}(r,\la)\wh \clp_2(\la)\big)
\\ & \label{B4!}
\times\big(\wh \clw_{11}(r,\la)\wh\clp_1(\la)+\wh\clw_{12}(r,\la)\wh \clp_2(\la)\big)^{-1},
\end{align}
where the pairs $\{\wh \clp_1,\wh \clp_2\}$ of the $p\times p$ matrix functions are nonsingular, with property-$J$, and
\begin{align} &       \label{B4!+}
\{\wh \clw_{ik}(r,\la) \}_{i.k=1}^2=\wh\clw(r,\la)=\wh W(r,\la)^{-1}.
 \end{align} 
In view of \eqref{B4}--\eqref{B4!+}, it is easy to see that the function $\wh \phi(r,\la)$,  generated by the pair 
$\{\wh \clp_1,\wh\clp_2\}$ (with property-$J$) of the form $\begin{bmatrix}\wh \clp_1(\la) \\ \wh \clp_2(\la)\end{bmatrix} =
\Th\begin{bmatrix}\clp_1(\la) \\ \clp_2(\la)\end{bmatrix}$, 
is connected with the function $\phi(r,\la)$ given by \eqref{P7}  by a simple linear-fractional transformation
\begin{align} &       \label{B5}
\wh \phi(r,\la)=\I \big(I_p+\phi(r,\la)\big)\big(I_p-\phi(r,\la)\big)^{-1}.
 \end{align} 
Since $\phi(r,\la)$ are contractions, the matrix functions  $\wh \phi(r,\la)\in\wh \cln(r)$ belong to Herglotz class, that is, $\I\big(\wh \phi(r,\la)^*-\wh \phi(r,\la)\big)\geq 0$.
Hence, the matrix functions  $\wh \phi(r,\la)$ admit Herglotz representation
\begin{align} & \label{B8}
\wh \phi(r,\la)=\mu \la+\nu+\int_{-\infty}^{\infty}\left(\frac{1}{t-\la}-\frac{t}{1+t^2}\right)d\tau(t), \quad \mu\geq 0, \quad \nu=\nu*,
\end{align}
where $\tau(t)$ is a $p\times p$ matrix function such that  $\tau(t_1)\geq \tau(t_2)$ for $t_1>t_2$ (i.e., $\tau$ is  monotonically increasing) and
\begin{align} & \label{B9}
\frac{d\tau(t)}{1+t^2}<\infty .
\end{align}
\begin{Dn}\label{DnHatW}
The matrix functions $\wh \phi(\la)\in\wh \cln(r)$, that is,  $\wh \phi(\la)$ of the form \eqref{B4!}  are called Weyl functions of the canonical system
\eqref{B1}, \eqref{B2} on $[0,r]$.
\end{Dn}

In case of the correspondence \eqref{B3}, the
operators $K=\I \wh \b(x)J\int_0^x\wh \b(t)^*\cdot dt$ and $V$ considered in this section coincide with the operators $K$ and $V$ from Section \ref{Prel} and Appendix~\ref{Simil}.
(Thus, we don't write $\wh K$ or $\wh V$.) Under condition 
\begin{align} &       \label{B4+}
\det \wh\b_2(0)\not=0,
 \end{align} 
the  operators $\wh E$ and $\wh V_0$ (acting in  $L_2^{p}(0, r)$) are introduced similar to $E$ and $V_0$, respectively:
\begin{align} & \label{B11}
 \wh E=V\wh V_0, \quad \wh V_0 f=\wh\b_2(0)f+\int_0^x\wh \clv_0(x-t)f(t)dt, \quad \wh \clv_0(x)=\big(V^{-1}\wh \b_2\big)^{\prime}(x).
\end{align} 
Compare  \eqref{B11} with \eqref{A36}. Note that the 
complete analog of Lemma~\ref{LaA2} is valid in our case and its proof coincides with the
proof of  Lemma~\ref{LaA2}.
Hence, an analogue of  Proposition \ref{PnM} is valid (see below).
\begin{Pn} \label{PnM'} The operator $\wh E$ given by \eqref{B11}
satisfies the equalities
\begin{align} & \label{B11+}
K=\wh EA \wh E^{-1}, \quad \wh E^{-1}\b_2 \equiv I_{p}.
\end{align}
Moreover, the operators $\wh E$ and $\wh E^{-1}$  map  $\clu^{p}[0,r]$ into $\clu^{p}[0,r]$.
\end{Pn}
Here, the proof of the second equality in \eqref{B11+} coincides
with the proof of the second equality in \eqref{A42}. 

Similar to the operators $S$ and $\Phi_1$, we introduce the operators $\wh S$ and  $\wh \Phi_1$ 
$\big(\wh \Phi_1\in \bB\big(\BC^{p},\, L_2^{p}(0, r)\big)$:
\begin{align} & \label{B10}
\wh S=\wh E^{-1}(\wh E^*)^{-1}, \quad \wh \Phi_1 h=\wh \Phi_1(x)h, 
\quad \wh \Phi_1(x)=\big(\wh E^{-1}\wh \b_1\big)(x).
\end{align}
Relations \eqref{B11+} and \eqref{B10} yield an analog
of  the operator identity \eqref{F3}:
\begin{align} & \label{F3'}
A\wh S-\wh SA^*=\I\wh \Pi J \wh \Pi^*, \quad \wh \Pi:=\begin{bmatrix}\wh \Phi_1 &  \Phi_2\end{bmatrix}.
\end{align}
The operator $\Phi_2$  in \eqref{F3'} coincides with the embedding $\Phi_2$ introduced in Section \ref{Fund}. 
The transfer matrix function $\wh w_A$ has the form
\begin{align} & \label{F6'}
\wh w_A(\la)=\wh w_A(r,\la)=I_{2p}-\I J \wh \Pi^*\wh S^{-1}(A-\la I)^{-1}\wh \Pi,
\end{align}
and the equality 
\begin{align} & \label{F14'}
\wh W(r,\la)=\wh w_A(r,1/\la)
\end{align}
follows similar to \eqref{F14} from the continuous
factorisation theorem \cite[p. 40]{SaL2} (or from its particular case \cite[Theorem 1.20]{SaSaR}).

{\bf 2.} Now, we are ready to study the high energy asymptotics of the Weyl functions $\wh \phi(\la)$.
Let $\wh \phi(\la)\in \wh \cln(r)$. Then, the formulas \eqref{B4!} and \eqref{B4!+} together with the property-$J$ of the pairs $\{\wh\clp_1,\wh \clp_2\}$
yield
\begin{align} & \label{B20}
\begin{bmatrix}I_{p} & \I\wh \phi(\la)^*\end{bmatrix}\wh W(r,\la)^*J \wh W(r,\la)\begin{bmatrix}I_{p} \\ -\I \wh \phi(\la)\end{bmatrix}\geq 0 .
\end{align}
Taking into account \eqref{F3'}--\eqref{F14'} and using again \cite[(1.88)]{SaSaR}, we obtain an analogue of \eqref{M1}:
\begin{align} & \label{M1'}
\wh W(r,\la)^*J\wh W(r,\la)=J+\I(\la -\ov{\la})\wh \Pi^*(I-\ov{\la}A^*)^{-1}\wh S^{-1}(I-{\la}A)^{-1} \wh \Pi.
\end{align}
It easily follows from \eqref{B20} and \eqref{M1'} that
\begin{align} & \label{B21}
\frac{\wh \phi(\la)-\wh\phi(\la)^*}{\la -\ov{\la}}\geq \begin{bmatrix}I_{p} &\I \wh \phi(\la)^*\end{bmatrix}\wh \Pi^*(I-\ov{\la}A^*)^{-1}\wh S^{-1}(I-{\la}A)^{-1} \wh \Pi
\begin{bmatrix}I_{p} \\ -\I \wh \phi(\la)\end{bmatrix}.
\end{align}
Note that $\mu=0$ in the Herglotz representation \eqref{B8} of $\wh \phi\in \wh\cln(r) $ (see \cite[(4.1)]{ALSinvpr}).
Hence, \eqref{B8} and \eqref{B9} imply that the norm of the left-hand side of \eqref{B21}
tends to zero in any angle $\de <\arg(\la)<\pi-\de$ $(\de>0)$ when $\la \to \infty$ in that angle.
Therefore, taking into account \eqref{M7} and \eqref{M8}, we similar to \eqref{M9} derive
\begin{align} & \label{B22}
\wh \phi(z^2)=-z\int_0^r\E^{\I z t}\wh \Phi_1(t)dt+ o\big(z\E^{\I z r}\big) \quad {\mathrm{for}} \quad |z|\to \infty,  
\end{align}
where $ \de/2<\arg(z)<(\pi-\de)/2$.
\begin{Tm}\label{TmHEAH}
Let canonical system \eqref{B1}, \eqref{B2} be given, such that $\wh\b(x)\in \clu^{p\times 2p}[0,r]$ and \eqref{B4+} holds.
Then, the Weyl functions $\wh \phi~\in~\wh \cln(r)$ admit $($for any $\de>0)$ the representation \eqref{B22}, where $\wh \Phi_1(x)\in \clu^{p\times p}[0,r]$ is given by  \eqref{B10}. The matrix function $\wh \Phi_1(x)$ is uniquely determined on $[0,r]$ by  the representation \eqref{B22}. 
\end{Tm}
\begin{proof}. The remaining proof that $\wh \Phi_1(x)$ is uniquely determined on $[0,r]$ by $\wh \phi(\la)$ is similar to the proof of Theorem \ref{TmHEA}.
Indeed, let another continuous matrix function $\breve \Phi(x)$ satisfy \eqref{B22}. Then, we have
\begin{align} & \label{M11'}
\om(z):=\int_0^r\E^{ \I z (t-r)}\big(\wh \Phi_1(t)-\br \Phi(t)\big)dt=o(1)  \quad {\mathrm{for}} \quad |z|\to \infty, 
\end{align}
where $ \de/2<\arg(z)<(\pi-\de)/2$. Clearly, $\om(z)$ is bounded in the lower half-plane $\Im(z)\leq 0$.
Hence, taking into account \eqref{M11'} and using  Phragmen--Lindel\"of theorem, we see that the entire function
$\om(z)$ is bounded in $\BC$ and tends to zero on some rays. Therefore, $\om(z)\equiv 0$, and so
$\wh \Phi_1(x)\equiv \breve \Phi(x)$ on $[0,r]$.
\end{proof}
Let us consider a simple example.
\begin{Ee} \label{Ee6.4} In the case $\wh S=I$, we have $A\wh S-\wh S A^*=\int_0^r(t-x)\, \cdot \, dt$.
Thus, the operator identity \eqref{F3'} holds for
\begin{align} & \label{E1}
\wh S=I, \quad \wh\Phi_1(x)=\I x I_p, \quad  \Phi_2(x)\equiv I_p.
\end{align}
Recall that $\wh \Pi$ in \eqref{F3'} is determined by the equality  $\wh \Pi \, h=\begin{bmatrix}\wh \Phi_1(x) & \Phi_2(x)\end{bmatrix} h$.
According  to the continuous factorisation theorem  \cite[Theorem 1.20]{SaSaR} $($and to \eqref{E1}$)$, $\wh W(x,\la)=\wh w_A(x,1/\la)$
is the normalised fundamental solution of the canonical system \eqref{B1} where
\begin{align} & \label{E2}
\wh H(x)=\wh\b(x)^*\wh\b(x), \quad \wh \b(x)=\begin{bmatrix}\wh \Phi_1(x) &  \Phi_2(x)\end{bmatrix}=\begin{bmatrix}\I x I_p & I_p\end{bmatrix}.
\end{align}
In view of \eqref{E2}, relations \eqref{B2} are satisfied. This example $($for different purposes$)$ was introduced in
\cite[p. 164]{SaL2}.
Taking into account the equalities \eqref{F6'}, \eqref{F14'} as well as \eqref{M4}, \eqref{M5}, \eqref{E1},
we calculate in a standard way the $p\times p$ blocks $\wh W_{ik}$ of $\wh W$:
\begin{align} & \nn
\wh W_{11}(r,\la)=\frac{1}{2}\big(\E^{\I\sqrt{\la} \, r }+\E^{-\I\sqrt{\la} \, r }\big)I_p, \quad  \wh W_{12}(r,\la)=\frac{\sqrt{\la}}{2}\big(\E^{\I\sqrt{\la} \, r }-\E^{-\I\sqrt{\la} \, r }\big)I_p, 
\\ & \nn
\wh W_{21}(r,\la)=\frac{1}{2\I}\Big(r\big(\E^{\I\sqrt{\la} \, r }+\E^{-\I\sqrt{\la} \, r }\big)
-\frac{1}{\I \sqrt{\la}}\big(\E^{\I\sqrt{\la} \, r }-\E^{-\I\sqrt{\la} \, r }\big)\Big)
I_p, 
\\ & \nn
\wh W_{22}(r,\la)=\frac{1}{2}\Big(\E^{\I\sqrt{\la} \, r }+\E^{-\I\sqrt{\la} \, r }-\I\sqrt{\la}\, r\big(\E^{\I\sqrt{\la} \, r }-\E^{-\I\sqrt{\la} \, r }\big)\Big)I_p.
\end{align}
By virtue of \eqref{B4!+} and \eqref{M1'} we have $\wh \clw(r,\la)=J\wh W(r,\ov{\la})^*J$, and so the expressions for
$\wh W_{ik}$ above imply the following representation of $\wh \clw$:
\begin{align} & \label{E3}
\wh\clw(r,\la)=\frac{1}{2}\E^{-\I z r }\begin{bmatrix}(\I z r+1) I_p & z I_p \\ \frac{1}{z}(\I z r+1) I_p  & I_p\end{bmatrix}+O(z\E^{\I z r}),
\end{align}
where $\la\in \BC_+$, $z=\sqrt{\la}$ and the branch  $\sqrt{\la}$ is chosen so that \\ $0<\arg(z)<\pi/2$. We choose a simple
nonsingular pair $($with property-$J):$
\begin{align} & \label{E4}
\wh \clp_1\equiv I_p, \quad \wh \clp_2\equiv I_p.
\end{align}
It easily follows from \eqref{B4!}, \eqref{E3} and \eqref{E4} that in the angle \\ $ \de/2<\arg(z)<(\pi-\de)/2$ $($for any $\de>0)$
we have
\begin{align} & \label{E5}
\wh \phi(z^2)=\frac{\I}{z}I_p+O\big(\E^{2\I z r}\big) \quad (|z|\to \infty).
\end{align}

On the other hand, substitution $\wh \Phi_1(t)=\I t I_p$ in the expression on the right-hand side of \eqref{B22} yields
\begin{align} & \label{E6}
-z\int_0^r\E^{\I z t}\wh \Phi_1(t)dt+ o\big(z\E^{\I z r}\big) =\frac{\I}{z}I_p+o\big(z\E^{\I z r}\big) \quad (|z|\to \infty).
\end{align}
Formula \eqref{B22} for our example is immediate from \eqref{E5} and \eqref{E6}.
That is, $\wh \Phi_1(x)=\I x I_p$ is the matrix function given by the last equality in \eqref{B10}.
\end{Ee}
Let us consider canonical system of the form \eqref{1.2}, \eqref{1.3} corresponding to system \eqref{B1}, \eqref{E2},
that is, canonical  system \eqref{1.2}, \eqref{1.3} where
\begin{align} & \label{E7}
m_1=m_2=p, \quad \b(x)=\wh \b(x)\Theta=\begin{bmatrix}\I x I_p & I_p\end{bmatrix}\Theta .
\end{align}
\begin{Ee}\label{Ee6.5}
Relations \eqref{B5} and  \eqref{E5} show that the Weyl function of the system \eqref{1.2}, \eqref{1.3}, \eqref{E7}
generated by the pair $\clp_1=I_p$, $\clp_2=0$ $($corresponding to the pair \eqref{E4}$)$ has the form
\begin{align} & \label{E8}
\phi(z^2)=\frac{1-z}{1+z}I_p+O\big(\E^{2\I zr}\big).
\end{align}
Thus, it is easily checked directly that
\begin{align} & \label{E9}
\Phi_1(t)=(2\E^{\I t}-1)I_p
\end{align}
satisfies \eqref{M9}. $($In other words,  we uniquely recovered $\Phi_1$ from $\phi$ in this example.$)$
\end{Ee}
\appendix
\section{Linear similarity transformation} \label{Simil}
\setcounter{equation}{0}
{\bf 1.} In this appendix, we consider linear operators $V$ and $V^{-1}$, which appear in the similarity transformation \eqref{AC3},
in greater detail. The operators $V$  may be constructed quite similar to the operators $V$ in \cite[Appendix C]{ALSstring}, where
some results of \cite{SaL0} are further developed.
First, we note that (similar to the case $m_1=m_2=p$ in \cite[Appendix C]{ALSstring}) the $m_2\times m_2$ integral kernel 
$\clv(x,\zeta)$ of the operator $V$ constructed in the way of  \cite[Appendix C]{ALSstring} admits representation
\begin{align} 
& \label{A1}
 \clv(x,\ze):=\sum_{k=1}^{\infty}\clv_k(x,\ze), 
\end{align}
where $\clv_1(x,\zeta)$ has the form
\begin{align}
 \nn
\clv_1(x,\zeta)=&\frac{1}{2}\left(\int_0^{(x+\zeta)/2}u_4(t)dt + \int_0^{(x-\zeta)/2}u_4(t)dt -\int_{(x+\zeta)/2}^x\Breve \clf(t, x-t+\zeta)dt\right.
\\  & \label{A2}
\left. 
- \int_{(x-\zeta)/2}^{x-\zeta}\Breve \clf(t, x-t-\zeta)dt-\int_{(x-\zeta)}^x\Breve \clf(t, \zeta+t-x)dt\right),
\end{align}
and the matrix functions $\clv_k(x,\zeta)$ (for $k>1$) have the form 
\begin{align} \nn
2\clv_k(x,\zeta)=&\int_{x-\zeta}^x\int_{\zeta+t-x}^tu_4(s)\clv_{k-1}(s,\zeta+t-x)dsdt
\\ \nn &
+\int_{(x+\zeta)/2}^x\int_{\zeta+x-t}^t u_4(s)\clv_{k-1}(s,\zeta+x-t)dsdt
\\  \label{A3} &
+\int_{(x-\zeta)/2}^{x-\zeta}\int_{x-t-\zeta}^t u_4(s)\clv_{k-1}(s,x-t-\zeta)dsdt
\\ \nn &
-\int_{x-\zeta}^{x}\int_{\zeta+t-x}^t \int_{\zeta+t-x}^s \clf(s,\eta)\clv_{k-1}(\eta,\zeta+t-x)d\eta dsdt
\\ \nn &
-\int_{(x+\zeta)/2}^{x}\int_{\zeta+x-t}^t \int_{\zeta+x-t}^s \clf(s,\eta)\clv_{k-1}(\eta,\zeta+x-t)d\eta dsdt
\\ \nn &
-\int_{(x-\zeta)/2}^{x-\zeta}\int_{x-t-\zeta}^t \int_{x-t-\zeta}^s \clf(s,\eta)\clv_{k-1}(\eta,x-t-\zeta)d\eta dsdt.
\end{align}
The $m_2\times m_2$ matrix functions $u_4$, $\clf$ and $\breve \clf$ in \eqref{A2} and \eqref{A3} satisfy relations
\begin{align}  \label{A4}&
u_4(t)=u_4(t)^*, \quad u_4\in L_2^{m_2\times m_2}(0,r); \quad \clf(s,\eta)=h_1(s)h_2(\eta), 
\\   \label{A5} &
h_1\in L_2^{m_2\times m}(0,r),
\quad h_2\in L_2^{m\times m_2}(0,r); \quad
\Breve \clf(t, \eta):=\int_{\eta}^t\clf(s,\eta)ds.
\end{align}
Here, the notations $h_k$ slightly differ from the corresponding notations in \cite{ALSstring}. For $C(r)=C>0$
such that
\begin{align} & \label{A6}
 \int_0^r \|h_k(t)\|dt \leq C \quad (k=1,2),  \quad \int_0^r \|u_4(t)\|dt \leq C^2,
\\ & \label{A7}
\sup_{0\leq \zeta \leq x \leq r} \|\clv_1(x,\zeta)\|\leq C,
\end{align}
the following inequalities are valid:
\begin{align} & \label{A8}
\|\clv_k(x,\zeta)\|\leq \frac{(3C^2)^{k-1}}{(k-1)!}Cx^{k-1} \quad (k\geq 1).
\end{align}
Similar to \cite[(C.22)]{ALSstring}, we have the equality
\begin{align} & \label{A8!}
u(0)=I_{m_2}.
\end{align}

{\bf 2.} Using the above-mentioned estimates and equalities, we prove our next proposition, which completes Theorem \ref{TmSim}.
\begin{Pn} \label{PnA1} Let the conditions of Theorem \ref{TmSim} hold. Then, the similarity transformation operators $V$ and $V^{-1}$ in \eqref{AC3}
$($constructed by the procedure from \cite[Appendix C]{ALSstring}$)$ map  vector  functions $f\in \clu^{m_2}[0,r]$ 
 into  $\clu^{m_2}[0,r]$ $($where $ \clu^{m_2}[0,r]= \clu^{m_2\times 1}[0,r]$ is introduced in \eqref{1.3+}$)$.
\end{Pn}
\begin{proof}.  Step 1.   Since $A \in\bB(L_2^{m_2}(0,r))$ is the operator of the squared integration multiplied by  $-1$, the relation $f\in \clu^{m_2}[0,r]$
is equivalent to the representation
\begin{align} & \label{A8+}
f(x)=-Af^{\prime\prime}+xf^{\prime}(0)+f(0).
\end{align}
In the first step, we show that $\Big(V\big(xf^{\prime}(0)\big)\Big)\in \clu^{m_2}[0,r]$, i.e.,
\begin{align} & \label{A8++}
\Big(V\big(xf^{\prime}(0)\big)\Big)^{\prime\prime}\in L_2^{m_2}(0,r).
\end{align}
For this purpose, we rewrite the first equality in \eqref{AC3} as $KV=VA$, where $V$ is the
integral operator given by the second equality in \eqref{AC3}, and
 $K$ is given in \eqref{P9}. Moreover, from the relations \eqref{A1}--\eqref{A8} it is easy to see that
 $\clv(x,\zeta)$ is continuous on the triangular $x \geq \zeta$.
Hence, in terms of the integral kernels, the equality 
$KV=VA$ (after we change the order of integration) is equivalent to
\begin{align}\nn
\I \b(x)j\big(\b(\ze)^*u(\ze)+\int_{\ze}^x\b(t)^*u(t)\clv(t,\ze)dt\big)=&(\ze -x)u(x)+u(x)\int_{\ze}^x(\ze - t)
\\ & \label{A9} \times
\clv(x,t)dt.
\end{align}

Sometimes, for convenience, we write  $Vg(t)$ (instead of $Vg(x)$).
Setting in \eqref{A9} $\ze=0$, we derive
\begin{align} & \label{A10}
V(tI_{m_2})=-\I \b(x)j\big(\b(0)^*+\int_{0}^x\b(t)^*u(t)\clv(t,0)dt\big),
\end{align}
where $V$ is applied to $tI_{m_2}$ columnwise.
We also note that according to \eqref{A2}, \eqref{A4} and \eqref{A5} we have
\begin{align} & \label{A11}
\clv_1(x,0)=\int_0^{x/2}u_4(t)dt-\int_0^{x/2}\left(\int_y^{x-y}h_1(s)ds\right)h_2(y)dy,
\\ &\label{A11+}
\frac{d}{dx}\clv_1(x,0)=\frac{1}{2}u_4\Big(\frac{x}{2}\Big)-\int_0^{x/2}h_1(x-y)h_2(y)dy\in L_2^{m_2\times m_2}(0,r).
\end{align}
From \eqref{A3} and \eqref{A4} we obtain the relations
\begin{align} \nn
\clv_k(x,0)=&\int_0^{x/2}\int_y^{x-y}u_4(s)\clv_{k-1}(s,y)dsdy
\\ & \label{A12}
-\int_0^{x/2}\int_y^{x-y}h_1(s)\int_y^s h_2(\eta)\clv_{k-1}(\eta,y)d\eta dsdy,
\\  \nn
\frac{d}{dx}\clv_k(x,0)=&
\int_0^{x/2}u_4(x-y)\clv_{k-1}(x-y,y)dy
\\ \label{A13} &
-\int_0^{x/2}h_1(x-y)\int_y^{x-y} h_2(\eta)\clv_{k-1}(\eta,y)d\eta dy,
\end{align}
for $k>1$. Formulas \eqref{A1}, \eqref{A8}, \eqref{A11+}, and \eqref{A13} show
that $\clv(x,0)$ is differentiable and
\begin{align} & \label{A14}
\frac{d}{dx}\clv(x,0)\in L_2^{m_2\times m_2}(0,r).
\end{align}
In view of \eqref{A10} and \eqref{A14}, $V(t I_{m_2})$  is two times differentiable
and 
\begin{align} & \label{A14+}
\frac{d^2}{dx^2}V(t I_{m_2})\in L_2^{m_2\times m_2}(0,r).
\end{align}
Thus, \eqref{A8++} holds.

Step 2. Next, we should show that $\clv(x,\ze)$ is differentiable with respect to $\ze$ on $[0,x]$.
Taking into account \eqref{A2}, \eqref{A4} and \eqref{A5}, it is easy to see that
\begin{align} & \label{A15}
\frac{\p}{\p \ze}\clv_1(x,\ze)=\frac{1}{4}u_4\left(\frac{x+\ze}{2}\right)-\frac{1}{4}u_4\left(\frac{x-\ze}{2}\right)\in L_2^{m_2\times m_2}(0,x).
\end{align}
It follows from \eqref{A3} and  \eqref{A4} that $\clv_2(x,\ze)$ is differentiable and for some $\wh C=\wh C(r)>0$ we  have
\begin{align} & \label{A16!}
\sup_{0\leq \ze\leq x\leq r}\left\|\frac{\p}{\p \ze}\clv_2(x,\ze)\right\|<\wh C.
\end{align} 
In view of \eqref{A2} and \eqref{A3}, for $\clv_k(x,x):=\lim_{\ze\to x-0}\clv(x,\ze)$ we also have
\begin{align} & \label{A16}
\clv_1(x,x)=\frac{1}{2}\int_0^xu_4(t)dt, \quad \clv_k(x,x)=0 \quad (k>1).
\end{align}
If $k>2$ and $\frac{\p}{\p \ze}\clv_{k-1}(x,\ze)$ exists and is bounded, simple calculations using \eqref{A3}, \eqref{A4} and \eqref{A16} show
that
\begin{align}  \nn
2\frac{\p}{\p \ze}\clv_k( x,\ze)=&\int_{x-\ze}^x\int_{\ze+t-x}^t u_4(s)\frac{\p}{\p \ze}\clv_{k-1}(s,\ze+t-x)dsdt
\\  \nn &
+\int_{(x+\ze)/2}^x\int_{\ze+x-t}^t u_4(s)\frac{\p}{\p \ze}\clv_{k-1}(s,\ze+x-t)dsdt
\\  \label{A17} &
+\int_{(x-\ze)/2}^{x-\ze}\int_{x-t-\ze}^t u_4(s)\frac{\p}{\p \ze}\clv_{k-1}(s,x-t-\ze)dsdt
\\  \nn &
-\int_{x-\ze}^{x}\int_{\ze+t-x}^t\int_{\ze+t-x}^s \clf(s,\eta)\frac{\p}{\p \ze}\clv_{k-1}(\eta,\ze+t-x)d\eta dsdt
\\  \nn &
-\int_{(x+\ze)/2}^{x}\int_{\ze+x-t}^t\int_{\ze+x-t}^s \clf(s,\eta)\frac{\p}{\p \ze}\clv_{k-1}(\eta,\ze+x-t)d\eta dsdt
\\  \nn &
-\int_{(x-\ze)/2}^{x-\ze}\int_{x-t-\ze}^t\int_{x-t-\ze}^s \clf(s,\eta)\frac{\p}{\p \ze}\clv_{k-1}(\eta,x-t-\ze)d\eta dsdt.
\end{align}

Now, we choose $\wh C$ in \eqref{A16!} such that $\wh C>C$, and so \eqref{A6} is valid for $\wh C$.
Let us show by induction that
\begin{align} & \label{A18}
\left\|\frac{\p}{\p \ze}\clv_k( x,\ze)\right\| \leq \frac{(3\wh C^2)^{k-2}}{(k-2)!}\wh Cx^{k-2} \quad (k\geq 2).
\end{align}
Indeed, \eqref{A18} holds for $k=2$. If \eqref{A18} is valid for $k-1$, we derive from \eqref{A4}, \eqref{A17}
and the inequality $\wh C>C$ that \eqref{A18} holds for $k$. 

Finally, relations \eqref{A1}, \eqref{A15} and \eqref{A18} imply the differentiability of $\clv(x,\zeta)$ with respect
to $\ze$ and the equality
\begin{align} & \label{A19}
\frac{\p}{\p \ze}\clv(x,\ze)=\sum_{k=1}^{\infty}\frac{\p}{\p \ze}\clv_k(x,\ze).
\end{align}
Moreover, the equalities \eqref{A15} and \eqref{A17} yield
\begin{align} & \label{A20}
\frac{\p}{\p \ze}\clv_k(x,\ze)\Big|_{\ze=0}\equiv 0 
\end{align}
for $k=1$ and for $k>2$. After slightly longer calculations, using \eqref{A3} one obtains $\frac{\p}{\p \ze}\clv_2(x,\ze)\Big|_{\ze=0}\equiv 0$.
Thus, \eqref{A20} holds for all $k>0$. Therefore, taking into account  
 \eqref{A19}, we derive
\begin{align} & \label{A21}
\frac{\p}{\p \ze}\clv(x,\ze)\Big|_{\ze=0}\equiv 0.
\end{align}

Step 3. 
It follows from \eqref{A1} and \eqref{A16} that
\begin{align} & \label{A22}
\clv(0,0)=0.
\end{align}
Taking the derivatives (with respect to $\ze$) of both parts of \eqref{A9} at $\ze=0$ and using \eqref{A21} and \eqref{A22}, we rewrite the result in the form
\begin{align} & \label{A23}
V(I_{m_2})=\I \b(x)j(\b^*u)^{\prime}(0), \quad {\mathrm{i.e.,}} \quad  V\big(f(0)\big)\in \clu^{m_2}[0,r].
\end{align}
Moreover, using again the equality $KV=VA$ we easily obtain
\begin{align} & \label{A24}
\Big(VAg\Big) \in \clu^{m_2}[0,r] \quad {\mathrm{for}}\quad g(x) \in L_2^{m_2}(0,r).
\end{align}
Finally, relations \eqref{A8+}, \eqref{A8++}, \eqref{A23} and \eqref{A24} yield
\begin{align} & \label{A25}
Vf \in \clu^{m_2}[0,r] \quad {\mathrm{for}} \quad f \in \clu^{m_2}[0,r].
\end{align}

Step 4.  Let us consider the operator $V^{-1}$. For this purpose, we rewrite the first equality in \eqref{AC3} as $V^{-1}K=AV^{-1}$.
In view of the first equality in \eqref{1.4} and the second equality in \eqref{P9}, we obtain
\begin{align} & \label{A27}
\big(Kg\big)(0)=0, \quad  \big(Kg\big)^{\prime}(0)=0 \quad {\mathrm{for}} \quad g\in L_2^{m_2}(0,r).
\end{align}
Thus, presenting the function $Kg$ in the form \eqref{A8+} and using \eqref{A27}, we derive $Kg=A(-Kg)^{\prime\prime}$
or, equivalently,
\begin{align} & \label{A28}
K=A\wh K, \quad \wh K:=I-\I\b^{\prime\prime}(x)j\int_0^x\b(t)^*\,\cdot \, dt.
\end{align}
Since $\wh K$ has a semi-separable kernel, it is invertible and $\wh K^{-1}\in \bB\big(L_2^{m_2}(0, r)\big)$
(see explicit expressions \cite[(C.5)--(C.7)]{ALSstring} for the inverse operator $\wh K^{-1}$). Hence, using \eqref{A28} we can rewrite
$V^{-1}K=AV^{-1}$ in the form
\begin{align} & \label{A29}
V^{-1}A=AV^{-1}\wh K^{-1}.
\end{align}

From \eqref{1.4} and \eqref{A10} one can see that
\begin{align} & \label{A30}
\big(V(xI_{m_2})\big)(0)=0, \quad \big(V(xI_{m_2})\big)^{\prime}(0)=I_{m_2}.
\end{align}
Hence, representation \eqref{A8+} of $V(xI_{m_2})$ has the form 
\begin{align} & \label{A31}
V(xI_{m_2})=xI_{m_2}-A\big(V(xI_{m_2})\big)^{\prime\prime}.
\end{align}
Multiplying both parts of \eqref{A31} by $V^{-1}$, we easily rewrite the result
in the form
\begin{align} & \label{A32}
V^{-1}(xI_{m_2})=xI_{m_2}+V^{-1}A\big(V(xI_{m_2})\big)^{\prime\prime}.
\end{align}
In the same way, the equality \eqref{A23} (in view of \eqref{1.4}, \eqref{A8!} and \eqref{A8+}) yields
\begin{align} & \label{A33}
V\, I_{m_2}=I_{m_2}+x\big(V\, I_{m_2}\big)^{\prime}(0)- A\big(V \, I_{m_2}\big)^{\prime\prime},
\\  & \label{A34}
V^{-1}\, I_{m_2}=I_{m_2}-V^{-1}\big(x\big(V\, I_{m_2}\big)^{\prime}(0)\big)+V^{-1}A\big(V \, I_{m_2}\big)^{\prime\prime}.
\end{align}
Clearly, the right-hand side of \eqref{A29} (and so the left-hand side as well)
maps functions from $L_2^{m_2}(0, r)$  into $\clu^{m_2}[0,r]$.
Hence, \eqref{A32} implies that $V^{-1}(xI_{m_2}) \in \clu^{m_2}[0,r]$. Therefore,
\eqref{A34} implies that the matrix function $V^{-1}\, I_{m_2}$ also belongs to $\clu^{m_2}[0,r]$. Summing up, we have
\begin{align} & \label{A35}
V^{-1}Af^{\prime\prime}, \,          V^{-1}\big(xf^{\prime}(0)\big), \,
V^{-1} f(0)
\in \clu^{m_2}[0, r]
\end{align}
for any $f \in \clu^{m_2}[0, r]  $. Thus, using \eqref{A8+} we see that
the statement of the proposition regarding operator $V^{-1}$ is valid.
\end{proof}

{\bf 3.} Assuming that the conditions of Proposition \ref{PnA1} hold,  
we introduce operators $E$ and $V_0$ by the equalities
\begin{align} & \label{A36}
E:=VV_0, \quad V_0f=\b_2(0)f+\int_0^x\clv_0(x-t) f(t) dt, \quad \clv_0:=\big(V^{-1}\b_2\big)^{\prime}, 
\end{align}
where $\b_2$ is the $m_2\times m_2$ block of $\b$ (see \eqref{F16}) and $E$ is another triangular operator
\begin{align} & \label{A36+}
Ef=u(x)\b_2(0)f+\int_0^x\cle(x,t)f(t)dt.
\end{align}
\begin{La} \label{LaA2} The operator $V_0$  is invertible and commutes with $A:$
\begin{align} & \label{A37}
V_0A=AV_0.
\end{align}
Moreover, the operators $V_0$ and $V_0^{-1}$  map  $\clu^{m_2}[0,r]$ into $\clu^{m_2}[0,r]$.
 \end{La}
\begin{proof}.
Taking into account \eqref{1.4}, it is easy to see that $\det \b_2(0)\not =0$. Indeed, if $\det \b_2(0) =0$ we have $ \b_2(0)^*g =0$ for some vector $g\not=0$. We also
have $\b_2(x)\b_2(x)^*=\b_1(x)\b_1(x)^*$ (which is immediate from $\b j\b^*\equiv 0$), and so $ \b_2(0)^*g =0$ implies $ \b_1(0)^*g =0$
and $\b(0)^*g=0$. Since $g\not=0$, the last equality contradicts the equality $\b(0)^{\prime}j\b(0)^*=\I I_{m_2}$. The invertibility of $V_0$
follows from the invertibility of $\b_2(0)$ and from the boundedness of $\clv_0(x)$ in the matrix norm.

The commutation property \eqref{A37} is immediate from the commutation property of the integral term  of $V_0$ and $A$.
In fact, there is a more general commutation property for triangular convolution operators. Indeed, for any  operator $\cla=\int_0^x a(x-t)I_{m_2} \cdot dt$
(where $a$ is a scalar function and is assumed, for convenience, to be bounded) 
we have
\begin{align} & \label{A38}
\int_0^x\clv_0(x-t)\int_0^t a(t-s)f(s)dsdt=\int_0^x\int_s^xa(t-s)\clv_0(x-t)dtf(s)ds, \\
& \label{A39}
 \int_0^xa(x-t)\int_0^t \clv_0(t-s)f(s)dsdt=\int_0^x\int_s^xa(x-t)\clv_0(t-s)dtf(s)ds.
\end{align}
Using the corresponding  change of variables, it is easy to show that the integral kernels on the right-hand sides of \eqref{A38} and \eqref{A39} coincide:
\begin{align} \nn
\int_s^xa(t-s)\clv_0(x-t)dt&=\int_0^{x-s}a(y)\clv_0(x-y-s)dy
\\ & \label{A40}
=\int_s^xa(x-t)\clv_0(t-s)dt.
\end{align}
Thus, the left-hand sides of \eqref{A38} and \eqref{A39} coincide as well (and \eqref{A37} is proved).

According to Proposition \ref{PnA1} and the last equality in\eqref{A36}, $\clv_0(x)$ is differentiable
and its derivative belongs to $L_2^{m_2\times m_2}(0,r)$.
Hence, the direct differentiation shows that $V_0$ maps  $ \clu^{m_2}[0,r]$ into the same class.

It is easy to see that
\begin{align} & \label{A41}
\big(V_0(xI_{m_2})\big)(0)=0, \quad \big(V_0(xI_{m_2})\big)^{\prime}(0)=\big(V_0 I_{m_2}\big)(0)=\b_2(0).
\end{align}
Finally, the fact that $V_0^{-1}$ maps $\clu^{m_2}[0,r]$ into the same class follows from \eqref{A37}
(more precisely from the equivalent equality $V_0^{-1}A=AV_0^{-1}$) and from \eqref{A41}
similar to the case of $V^{-1}$ in the proof of Proposition \ref{PnA1}.
\end{proof}
\begin{Rk} \label{RkLa} The proof  of  Lemma \ref{LaA2} shows that a wide class of operators commutes with $A$,
and so the transformation operators $\wh V$ such that $K=\wh V A\wh V^{-1}$ are not uniquely defined.
In this context, the {\it normalised} transformation operator $E$  given by \eqref{A36} is important.
\end{Rk}
\begin{Pn} \label{PnM} The operator $E$ given by \eqref{A36}
satisfies the equalities
\begin{align} & \label{A42}
K=EAE^{-1}, \quad E^{-1}\b_2 \equiv I_{m_2}.
\end{align}
Moreover, the operators $E$ and $E^{-1}$  map  $\clu^{m_2}[0,r]$ into $\clu^{m_2}[0,r]$.
\end{Pn}
\begin{proof}. The first equality in \eqref{A42} is immediate from the definition of $E$ and formulas \eqref{AC3}  and \eqref{A37}.
In order to derive the second equality in \eqref{A42}, we note that the representation of $V$ in \eqref{AC3} and
the inequality \eqref{AC4} yield that
\begin{align} & \label{A43}
\big(V^{-1}f\big)(x)=u(x)^*f(x)+\int_0^x\clq(x,\xi)f(\xi)d\xi, \quad \sup_{0\leq \xi\leq x\leq r}\|\clq(x,\xi)\|<\infty .
\end{align}
In particular, in view of \eqref{A8!} we have 
\begin{align} & \label{A44}
\big(V^{-1}\b_2\big)(0)=\b_2(0).
\end{align}
Taking into account \eqref{A44} and the last two equalities in \eqref{A36}, we derive
\begin{align} & \label{A45}
\big(V_0 I_{m_2}\big)(x)=\b_2(0)+\int_0^x\clv_0(\xi)d\xi=\big(V^{-1}\b_2\big)(x),
\end{align}
and the second equality in \eqref{A42} follows.

Finally, the fact that the operators $E$ and $E^{-1}$  map   $\clu^{m_2}[0,r]$ into $\clu^{m_2}[0,r]$ follows from the corresponding properties of $V^{\pm1}$
and $V_0^{\pm1}$ (see Proposition \ref{PnA1} and Lemma \ref{LaA2}).
\end{proof}
\begin{Rk}\label{RkALT} It follows from the Remark \ref{RkLT} and formulas \eqref{A36} and \eqref{A36+}
that the integral kernel
$\cle(x,t)$ $($of $E)$ in the domain $0\leq t\leq x\leq \ell<r$ is uniquely determined by
$\b(x)$ on $[0,\ell]$ $($and does not depend on the choice of $\b(x)$ for $\ell<x<r$ and the choice of $r\geq \ell)$.
\end{Rk}


\begin{flushright}
Alexander Sakhnovich \\
Faculty of Mathematics,
University
of
Vienna, \\
Oskar-Morgenstern-Platz 1, A-1090 Vienna,
Austria, \\
e-mail: {\tt oleksandr.sakhnovych@univie.ac.at}

\end{flushright}


\begin{thebibliography}{AGKS}

\bibitem{ArD}
D.Z. Arov and H. Dym, {\it Bitangential direct and inverse problems for systems of integral and differential equations.} Cambridge University Press, Cambridge, 2012.

\bibitem{Atk}
F.V. Atkinson, On the location of the Weyl Circles. {\it Proc. Roy. Soc. Edinburgh} {\bf 88 A}
(1981), 345--356.

\bibitem{ClGe}
S. Clark and F. Gesztesy,  Weyl-Titchmarsh M-function asymptotics, local uniqueness results, trace formulas, and Borg-type theorems for Dirac operators. 
{\it Trans. Amer. Math. Soc.} {\bf 354} (2002), no. 9, 3475--3534.

\bibitem{ClGe2}
S. Clark and F. Gesztesy, On self-adjoint and J-self-adjoint Dirac-type operators: a case study. Recent advances in differential equations and mathematical physics.
In {\it Contemp. Math.} {\bf 412}. Amer. Math. Soc., Providence, RI, 2006, 103--140.



\bibitem{dBr}
L.  de Branges,  {\it Hilbert spaces of entire functions.} Prentice-Hall,  Englewood Cliffs, N.J., 1968.

\bibitem{DuK}
L.N. Dudko and I.I. Kal'mushevskii, {Conditions for similarity to the operator $I^2$ in terms of the characteristic matrix-valued function}. In Russian.
{\it Izv. Vyssh. Uchebn. Zaved. Mat.} (1976), no. 4,
38--46.






\bibitem{EKT}
J. Eckhardt, A. Kostenko, and G. Teschl, 
{Spectral asymptotics for canonical systems.} {\it J. Reine Angew. Math.} {\bf 736} (2018), 285--315.

\bibitem{Ev}
W.N. Everitt, On a property of the m-coefficient of a second-order linear differential equation. {\it J. London Math. Soc.} {\bf 4} (1972), 443--457.

\bibitem{GeB}
F. Gesztesy,
Inverse spectral theory as influenced by Barry Simon. Spectral theory and mathematical physics: a Festschrift in honor of Barry Simon's 60th birthday.
In {\it Proc. Sympos. Pure Math.} {\bf 76}, Part 2. Amer. Math. Soc., Providence, RI, 2007,  741--820.

\bibitem{GeS}
F. Gesztesy and A. Sakhnovich,  {The inverse approach to Dirac-type systems based on the A-function concept.} {\it J. Funct. Anal.} {\bf 279} (2020), Art. 108609.

\bibitem{GeSi}
F. Gesztesy  and  B. Simon, 
A new approach to inverse spectral theory.  II. {General real potentials and the connection to the spectral measure}. \textit{Ann. of Math. (2)} {\bf 152} (2000), 593--643.

\bibitem{GeSi2}
F. Gesztesy  and  B. Simon,
On local Borg-Marchenko uniqueness results. 
{\it Comm. Math. Phys.}  {\bf 211} (2000),  273--287. 

\bibitem{GoKr}
I.C. Gohberg  and  M.G. Krein,
\textit{Theory  and  applications  of Volterra operators  in  Hilbert  space.} Transl.  math.  monographs {\bf 24}. Amer. Math. Soc., Providence,  RI, 1970.

 \bibitem{GolMi}
 L. Golinskii and I. Mikhailova (edited by V.P. Potapov), 
 {Hilbert spaces of entire functions as a $J$-theory subject}.  In 
  {\it Oper. Theory Adv. Appl.} \textbf{95}.
 Birkh\"auser, Basel,   1997,  205--251.


\bibitem{KKL} 
I. Koltracht, B.   Kon, and L.  Lerer,
{Inversion of structured operators.} \textit{ Integral Equations Operator Theory} {\bf  20} (1994), 410--448.

\bibitem{Krein0}
M.G. Krein,  On the transfer function of a one-dimensional boundary problem
of the second order. In Russian. {\it Doklady Akad. Nauk SSSR (N.S.)} {\bf 88} (1953),
405--408.

\bibitem{Krein} 
M.G. Krein,  {Continuous analogues of propositions on polynomials orthogonal on the unit circle.}  In Russian.
{\it Dokl. Akad. Nauk SSSR (N.S.)} \textbf{105} (1955), 637--640. 


\bibitem{KLang} 
M.G. Krein and H. Langer, 
Continuation of Hermitian positive definite functions and related questions. 
{\it Integral Equations Operator Theory} {\bf 78} (2014), no. 1, 1--69. 

\bibitem{Langer}
H. Langer,  Transfer functions and local spectral uniqueness for Sturm-Liouville operators, canonical systems and strings. {\it Integral Equations Operator Theory} {\bf 85} 
(2016),  1--23. 

\bibitem{LW}
M. Langer and H. Woracek,  A local inverse spectral theorem for Hamiltonian systems. {\it Inverse Problems} {\bf 27} (2011), no. 5, Art.  055002.




\bibitem{Mog}
V. Mogilevskii, 
Spectral and pseudospectral functions of Hamiltonian systems: development of the results by Arov--Dym and Sakhnovich.
{\it Methods Funct. Anal. Topology} {\bf 21} (2015), no. 4, 370--402. 


\bibitem{Rem}
C. Remling, 
{\it Spectral theory of canonical systems.} 
De Gruyter, Berlin, 2018.
 

\bibitem{Rom}
R. Romanov,  {Order problem for canonical systems and a conjecture of Valent.} {\it Trans. Amer. Math. Soc.} {\bf 369} (2017),  1061--1078. 

\bibitem{RW}
R. Romanov and H. Woracek,  {Canonical systems with discrete spectrum.} {\it J. Funct. Anal.} {\bf 278} (2020), Art. 108318.
 
 \bibitem{Rov}
 J. Rovnyak and L.A. Sakhnovich, {Pseudospectral functions for canonical differential systems. II.} 
 In {\it Oper. Theory Adv. Appl.} {\bf 218}. Birkh\"auser/Springer, Basel, 2012, 583--612. 
 
\bibitem{Ryb}
A. Rybkin, 
On the trace approach to the inverse scattering problem in dimension one. 
{\it SIAM J. Math. Anal.} {\bf 32} (2001), 1248--1264. 

 \bibitem{ALS88}
A.L. Sakhnovich,  Asymptotics of spectral functions of an S-colligation. {\it Soviet Math. (Iz. VUZ)} {\bf 32} (1988), no. 9, 92--105. 

\bibitem{ALS90}
A.L. Sakhnovich,  A nonlinear Schr\"odinger equation on the semi-axis and a related inverse problem.
{\it Ukrainian Math. J.} {\bf 42} (1990), no. 3, 316--323.

 \bibitem{ALS02}
A.L. Sakhnovich, {Dirac type and canonical  systems: spectral and
Weyl-Titchmarsh functions, direct and inverse  problems}. {\it Inverse
Problems} {\bf 18} (2002), 331--348.


\bibitem{ALS15}
A.L. Sakhnovich,  Inverse problem for Dirac systems with locally square-summable potentials 
and rectangular Weyl functions. {\it J. Spectr. Theory} {\bf 5} (2015), 547--569.

\bibitem{ALS17}
A.L. Sakhnovich, {Dynamical canonical systems and their explicit solutions.} {\it Discrete Contin. Dyn. Syst.} {\bf 37} (2017),  1679--1689. 











\bibitem{ALSstring} 
A.L. Sakhnovich, {On the class of canonical systems corresponding to matrix string equations: general-type and explicit fundamental solutions and Weyl--Titchmarsh theory}.
{\it Doc. Math.} {\bf 26} (2021), 583--615.

\bibitem{ALSinvpr} 
A.L. Sakhnovich, On the solution of the inverse problem for a class of canonical systems corresponding to matrix string equations.
ArXiv:2104.07612




\bibitem{SaSaR}
A.L. Sakhnovich, L.A. Sakhnovich,  and I.Ya. Roitberg,   \textit{Inverse Problems and Nonlinear Evolution Equations. 
 Solutions, Darboux Matrices and Weyl--Titchmarsh Functions}. De Gruyter,  Berlin, 2013.
 
\bibitem{SaL0} 
L.A.  Sakhnovich, The spectral analysis of Volterra operators and some inverse problems. In Russian. {\it Dokl. Akad. Nauk SSSR (N.S.)} {\bf 115} (1957), 666--669. 

\bibitem{SaL1}
L.A. Sakhnovich,  {On  the  factorization  of  the  transfer
matrix function.} {\it Sov. Math. Dokl.} { \bf 17} (1976), 203--207.

\bibitem{SaL1+}
L.A. Sakhnovich,
{Equations with a difference kernel on a finite interval}.
{\it Russian Math. Surveys}  {\bf 35} (1980), 81--152.


\bibitem{SaL2-}
L.A. Sakhnovich, 
{Factorization  problems  and  operator identities.} 
{\it Russian Math. Surveys}   {\bf 41} (1986), 1--64.

\bibitem{LA94}
L.A. Sakhnovich, {The method of operator identities and problems in analysis.}  {\it St. Petersburg Math. J.} {\bf 5} (1994),  1--69.


\bibitem{SaL2}
L.A. Sakhnovich,  {\it Spectral theory of canonical differential
systems, method of operator identities.}  Birkh\"auser, Basel, 1999.

\bibitem{SaL15}
{L.A. Sakhnovich},  {\it Integral equations with difference kernels on finite intervals}. Second edition, revised and extended.
Birkh\"auser/Springer, Cham, 2015.
 
\bibitem{Si0} 
B. Simon.       
A new approach to inverse spectral theory. I. Fundamental formalism. \textit{Ann. of Math.},   {\bf 150} (1999), 1029--1057.
 


\bibitem{Su}
M. Suzuki, {An inverse problem for a class of canonical systems having Hamiltonians of determinant one}. {\it J. Funct. Anal.} {\bf 279} (2020), Art. 108699.

\bibitem{Wor}
H. Woracek, {Asymptotics of eigenvalues for a class of singular Krein strings}. {\it Collect. Math.} {\bf 66} (2015), 469--479.



\end{thebibliography}
\end{document}